\documentclass[11pt, a4paper]{article}
\usepackage{amsmath,amssymb,amsfonts}
\usepackage{verbatim}
\usepackage{a4wide}
\usepackage{psfig,epsfig,epic,eepic,graphicx}
\newtheorem{theorem}{Theorem}
\newtheorem{lemma}[theorem]{Lemma}
\newtheorem{proposition}[theorem]{Proposition}

\newtheorem{corollary}{Corollary}

\makeatletter

\@addtoreset{equation}{section}
\makeatother

\newcommand{\E}{{\mathbb E}}
\newcommand{\T}{{\mathbb T}}
\newcommand{\N}{{\mathbb N}}
\newcommand{\Z}{{\mathbb Z}}

\newcommand{\R}{{\mathbb R}}

\newcommand{\pa}{{\partial}}
\newcommand{\na}{{\nabla}}

\newcommand{\eps}{{\varepsilon}}

\def\div{\hbox{div  }}

\title{The Navier wall law at a boundary with random roughness} 
\author{David G\'erard-Varet 
\footnote{DMA/CNRS, Ecole Normale Sup\'erieure, 45 rue d'Ulm,75005
   Paris, FRANCE}}
\date{}

\begin{document}
\maketitle

\begin{abstract}
We consider  the Navier-Stokes equation in a  domain with
irregular boundaries. The irregularity is modeled by a spatially
homogeneous random process, 
with typical size $\eps \ll 1$. In the parent paper \cite{Basson:2007},
 we derived a homogenized 
boundary condition of Navier type as $\eps \rightarrow 0$. We show here
that for a large class of boundaries,  this Navier condition provides a
$O(\eps^{3/2} |\ln \eps|^{1/2})$ approximation  in $L^2$, instead of
$O(\eps^{3/2})$ for periodic irregularities. Our result
relies on the study of an auxiliary boundary layer system. Decay
properties of this boundary layer are deduced from a central limit theorem
for dependent variables. 
\end{abstract}

{{\em Keywords:} \small Wall laws, rough boundaries, stochastic 
homogenization, decay of correlations} 
 
\section{Introduction}
The concern of this paper is the effect of a rough boundary on a viscous
 fluid. In most situations of physical relevance, such effect can not be
 described in detail: either the precise shape of the roughness is unknown,
 or its spatial variations are too small for computational grids.
 Therefore, one may only hope to account for  the averaged effect
of the irregularities. This is the purpose of {\em wall laws}: the
irregular boundary is replaced by an artificial smoothed one, and an
artificial  boundary condition (a wall law) is prescribed there, that
 should reflect the mean impact of the roughness. 

\medskip
This paper is a  mathematical study of  wall laws, in the following
simple setting:
we consider a two-dimensional rough channel 
$$\Omega^\eps \: =  \: \Omega \cup \Sigma \cup R^\eps $$
where $\Omega = \R \times (0,1)$ is the {\em smooth part}, $R^\eps$ is the
rough part, and $\Sigma = \R \times \{0\}$ their interface. We assume that
the rough part has typical size $\eps$, that is 
$$ R^\eps \: = \: \left\{ x, \: x_2 >
\eps\omega\left(\frac{x_1}{\eps}\right) \right\} $$ 
for a  $K-$Lipschitz function $\omega : \R \mapsto (-1,0)$, $K > 0$. 
More will be assumed on the boundary function $\omega$ hereafter 
(see figure for an example of such a rough domain).

\begin{figure}
\begin{center}
\includegraphics[height = 5.5cm, width=8.5cm]{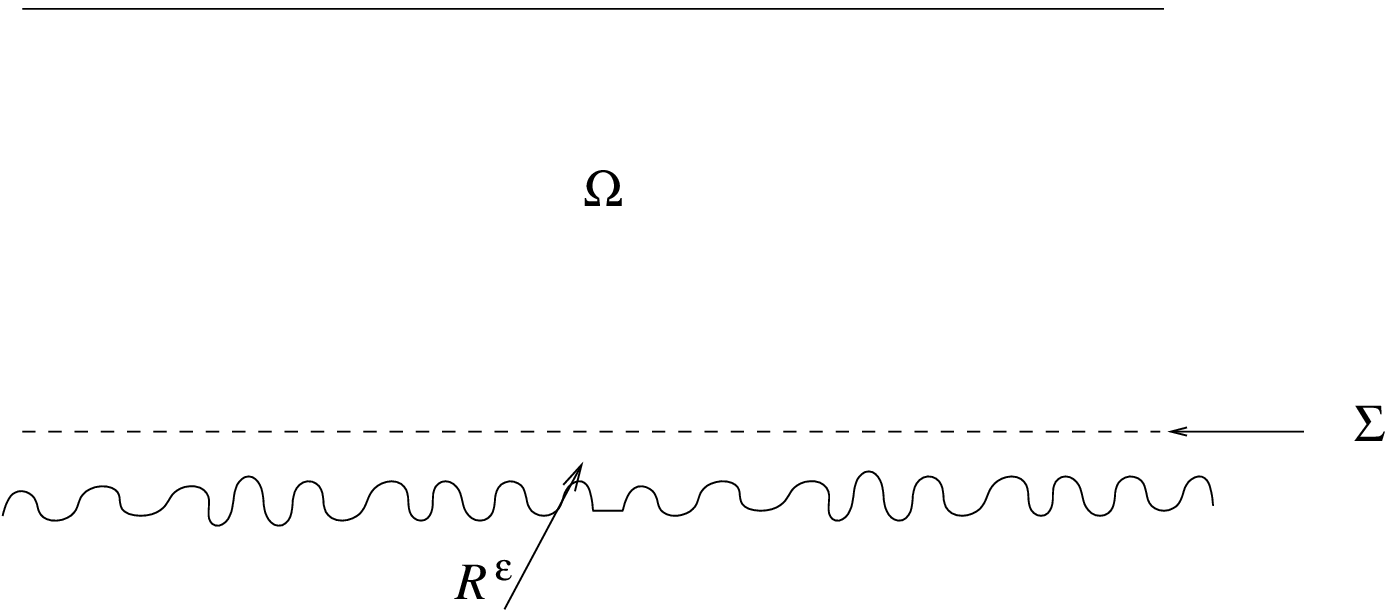}
\end{center}
\caption{The rough domain $\Omega^\eps$.}
\end{figure}

\medskip
 We assume that in this channel domain, 
the viscous fluid obeys to the stationary incompressible Navier-Stokes
equations: 
\begin{equation} \label{NSeps}
\left\{
\begin{aligned}
& - \Delta u + u \cdot \na u + \na p  = 0, \: x \in \Omega^\eps, \\
& \div u  = 0, \: x \in \Omega^\eps, \\
& \int_{\sigma^\eps} u_1   = \phi, \\
& u\vert_{\pa \Omega}  = 0, 
\end{aligned}
\right.
\end{equation}
where $\sigma^\eps$ denotes any vertical cross-section of $\Omega^\eps$
and $\phi > 0$. The third equation in \eqref{NSeps} expresses that a flux
$\phi$ is imposed across the channel. Note that this flux does not depend
on the cross-section, due to the incompressibility and  no-slip condition
at the boundary. We also stress that, up to minor changes, 
 we could apply our analysis to many variants of this problem, 
 notably to elliptic type systems or to unstationary  Navier-Stokes.

\medskip
In this simple setting, the search for  wall laws resumes to the following
problem: to find a boundary operator $B^\eps(x,D_x)$,
regular in $\eps$, acting  at the {\em artificial} boundary $\Sigma$,
such that the solution of  
 \begin{equation} \label{NS}
\left\{
\begin{aligned}
& - \Delta u + u \cdot \na u + \na p  = 0, \: x \in \Omega, \\
& \div u  = 0, \: x \in \Omega, \\
& \int_{\sigma} u_1   = \phi, \quad u\vert_{x_2=1}  = 0, \\
& B^\eps(x,D_x)u\vert_{\Sigma}  = 0
\end{aligned}
\right.
\end{equation}
approximates well the solution $u^\eps$ of \eqref{NSeps} in $\Omega$.

\medskip
This type of homogenization problems  has been considered in many
 mathematical works.  
 On wall laws for scalar elliptic equations, we refer to
 \cite{Achdou:1995}. On wall
 laws for fluid flows, see \cite{Achdou:1998a,
 Achdou:1998,Amirat:2001a,Jager:2001,Jager:2003, Bresch:2006}.
 See also \cite{Jager:2000} on porous
 boundaries. These works go along with more formal computations, grounded by
 empirical arguments ({\it cf} for instance \cite{Bechert:1989,
  Luchini:1995}).  
We finally  mention \cite{Gerard-Varet:2003b,Bresch:2005} for the study of
 roughness-induced effects  on geophysical systems.
  
All these studies  have been  carried under two assumptions:
\begin{itemize}
\item compact domains, for instance bounded channels or periodic in the
   variable $x_1$.
\item {\em periodic irregularities},
meaning that the boundary function $\omega$ is 
periodic.
\end{itemize}
The first restriction is just a  small mathematical convenience, that gives
 direct compactness properties through Rellich type theorems. The second
 assumption is of course a big simplification, both from the point of
 view of mathematics and physics. These assumptions were considerably
 relaxed in the recent article  \cite{Basson:2007}
 by A. Basson and the author. As the
 present note extends this article, we now describe shortly its main
 results and underlying difficulties.  
  
\medskip
In all  papers on wall laws, the starting point is a formal expansion of
 $u^\eps$:
$$ u^\eps(x) \: \sim \:  u^0(x) +  6 \phi \eps v(x/\eps) + \dots $$
Formally, the leading term  $u^0$ satisfies \eqref{NS} with the simple
no-slip condition 
\begin{equation} \label{Dirichlet} 
B^\eps(x,Dx) u \: :=  \: u = 0  \:\mbox{ at }\: \Sigma 
\end{equation}
The solution of this approximate system  is the famous Poiseuille flow :
\begin{equation*}
u^0(x) \: = \: \left(U(x_2),0\right), \quad U(x_2) \: = \:  6 \phi x_2
(1-x_2) 
\end{equation*}
Note that $u^0$ is defined in all $\R^2$. 
This zeroth order asymptotics can be mathematically justified, 
at least for small fluxes $\phi$: we prove in article \cite{Basson:2007}
\begin{theorem} \label{Direstimates}
There exists $\phi_0, \: \eps_0 >0$, such that for  all $\phi < \phi_0, \:
 \eps < \eps_0$, system \eqref{NSeps} has a unique solution $u^\eps$ in
$H^1_{uloc}(\Omega^\eps)$. Moreover, 
$$ \| u^\eps - u^0 \|_{H^1_{uloc}(\Omega^\eps)} \le C \sqrt{\eps}, \quad \|
u^\eps - u^0 \|_{L^2_{uloc}(\Omega)} \le  C'  \eps. $$
\end{theorem}
We stress that these estimates hold without further assumption on the
boundary: we only assume that $\omega$ has values in $(-1,0)$ and is
$K-$Lip. A look at the proof shows that the constants $C$ and $C'$ are  only
increasing functions of $K$.  

\medskip
Theorem \ref{Direstimates} expresses that the wall law
\eqref{Dirichlet} provides a $O(\eps)$ approximation of $u^\eps$ in
$L^2_{uloc}(\Omega)$. See also \cite{Jager:2001}
 for a similar result in a bounded channel. However, this wall law does not 
account for the behaviour of
$u^\eps$ near the boundary, and can therefore be refined. Indeed, as the
Poiseuille flow $u^0$ does not vanish at the lower part of  $\pa \Omega^\eps$, 
a boundary layer corrector $\eps \phi v(x/\eps)$ must be added to the
expansion. The (normalized) boundary layer $v = v(y)$ is defined on the 
rescaled infinite domain 
$$ \Omega^{bl} \: = \: \{ y, \: y_2 > \omega(y_1)  \} $$
and formally satisfies the following Stokes problem
\begin{equation} \label{BL}
\left\{
\begin{aligned}
& -\Delta v + \na q = 0, \: x \in \Omega^{bl}, \\
& \div v = 0, \: x \in \Omega^{bl},\\
& v(y_1, \omega(y_1)) \: = \: - (\omega(y_1),0).
\end{aligned}
\right.
\end{equation}
Note the inhomogeneous Dirichlet condition, that  cancels the trace of $u^0$.  

\medskip
{\em Although linear, the boundary layer system \eqref{BL} is
quite challenging}. First, the well-posedness is not
  clear. As the boundary function $\omega$ is not decreasing at infinity,
  one can expect only local integrability of the solution  $v$ in
  variable $y_1$. The derivation of local bounds is not
  obvious: the Stokes operator being vectorial, one can not use scalar tools
  such as the maximum principle or Harnack inequality. Moreover, as
  $\Omega^{bl}$ is unbounded in all directions, the
  Poincar\'e inequality (which allows to get $H^1_{uloc}$ estimates in the
  channel) is not available. Besides the well-posedness issue, the 
qualitative  properties  of $v$ seem also out of reach without further 
hypothesis.  

Under an assumption of periodic irregularities, the analysis of \eqref{BL}
becomes straightforward.
 If $\omega$ is say $L$ periodic in $y_1$, it is easy to show
 well-posedness in the space  
$$   \left\{ v \in H^1_{loc}(\Omega_{bl}), \: v \:\:  
L-\mbox{periodic in } y_1,
\: \int_0^L \int_{\omega(y_1)}^{+\infty} | \na v |^2 dy_2 dy_1 < 
+\infty \right\}. $$ 
Moreover, a simple Fourier transform in $y_1$ shows that 
$$ \| v(y) - v^\infty \| \le C \, e^{-\delta y_2/L}, \quad  v^\infty = (\alpha,
0), \quad  \alpha= \frac{1}{L} \int_0^L v_1(s) ds,  \quad \delta > 0,$$
that is exponential convergence to a constant field $v^\infty = (\alpha,
0)$ at infinity. 

\medskip
The constant $\alpha$ at infinity is then responsible  for a $O(\eps)$
tangential slip. Namely, chosing as a wall law  the Navier-slip condition 
\begin{equation} \label{Navier}
 B^\eps(x, D_x)v \:  = \:  \left( v_1 - \eps \alpha \pa_2 v_1, \:
  v_2  \right) = 0 \: \mbox{ at } \Sigma,
\end{equation}  
it can be shown (in this periodic framework) that the solution of \eqref{NS}
provides a $O(\eps^{3/2})$ approximation of $u^\eps$ in $L^2$. We refer to
\cite{Jager:2001} for all necessary details.
 The error estimate $\eps^{3/2}$ comes
 from the fact that the boundary layer term
satisfies $\| \eps (v(x/\eps)- (\alpha,0)) \|_{L^2} = O(\eps^{3/2})$.  

\medskip
 The periodicity hypothesis is a stringent one, and one may wonder if the
  use of  Navier slip condition can be justified in more general
 configurations.  This issue has been adressed rigorously in the recent
 article \cite{Basson:2007}.
  Inspired by the probabilistic modeling of heterogeneous
  media (see for instance \cite{Jikov:1994}), we
 considered irregularities that are not distributed periodically, but
 randomly, following a stationary stochastic process. Namely, the rough
  boundary is seen as a realization of a stationary spatial
  process. Following the well-known construction  of Kolmogorov, this
  amounts to consider the space 
$$ P \: = \: \left\{ \omega : \R \mapsto (-1,0), \: \omega \: K-\mbox{Lip}
  \right\} $$
of all possible rough boundaries, together with the cylindrical $\sigma-$
  field ${\cal C}$ (that is generated by the coordinates $\omega \mapsto
  \omega(t)$) 
  and with a stationary  measure
  $\pi$. Stationary means that $\pi$ is invariant by the group of translation 
$$ \tau_h : P \mapsto P, \quad \omega \mapsto \omega(\cdot + h). $$
As a consequence of this modeling, the domains $\Omega^\eps$, $\Omega^{bl}$,
as well as the velocity fields $u^\eps$ or  $v$  depend on the 
parameter $\omega$. As discussed earlier,
the existence result and estimates of theorem \ref{Direstimates} are
uniform on $P$. Moreover, it was  shown in article \cite{Basson:2007}
 that the function
$\omega \mapsto u^\eps(\omega, \cdot)$ (extended by $0$ outside
 $\Omega^\eps(\omega)$) is measurable as a function from $P$
to $H^1_{loc}(\R^2)$. 

\medskip
Using this probabilistic structure, we have been able to extend
partially the results of the periodic case. Key elements of our
analysis are:
\begin{itemize}
\item  the well-posedness of the boundary layer system, obtained in
 functional spaces encoding the relation
$$ v(\tau_h(\omega),y_1,y_2) \: = \: v(\omega, y_1+h, y_2). $$
\item the convergence of $v(\omega, y)$  to $(\alpha(\omega), 0)$
 as $y_2 \rightarrow + \infty$, both in $L^2(P)$ and  almost surely,
 locally uniformly in $y_1$. Such convergence is deduced from the ergodic
 theorem.   
\end{itemize}
More on the boundary layer system will be provided in the next sections. 
 As regards the Navier wall law \eqref{Navier}, the
main result of \cite{Basson:2007} resumes to  
\begin{theorem} \label{Navestimates}
There exists $\alpha = \alpha(\omega) \in L^2(P)$ such that the solution
$u^N$ of \eqref{NS}, \eqref{Navier} satisfies 
$$   \| u^\eps - u^N \|_{L^2_{uloc}(P \times \Omega)} = o(\eps). $$
\end{theorem}
We remind that
$ \| w \|_{L^2_{uloc}(P \times \Omega)} \: := \: sup_{x} \left( \int_P
\int_{B(x,1) \cap \Omega} | w |^2 dx dP \right)^{1/2}.$

\medskip
Theorem \ref{Navestimates} shows that a slip
condition of Navier type improves the approximation of $u^\eps$. As in the
periodic case, the random variable $\alpha$ in \eqref{Navier} comes from 
the convergence of the
boundary layer $v$. If the measure $\pi$ is
ergodic, $\alpha$ does not depend on $\omega$, as pointed out in
\cite{Basson:2007}.   

A natural concern about this result is the $o(\eps)$ bound, which is only
a slight improvement of the $O(\eps)$ in theorem
\ref{Direstimates}. A look at article \cite{Basson:2007}
 shows that this  poor bound is
due to the lack of information on the way $v$ converges at infinity. 
 Contrary to the periodic case, where convergence at exponential
rate is established, the simple use of the ergodic theorem does not yield
any speed rate. 

\medskip
The present paper aims at clarifying this point. Losely, {\em we will show
  that for a large 
class of boundaries, the Navier wall law provides a $O(\eps^{3/2}
|\ln(\eps)|^{1/2})$
 approximation of the real solution.} Namely, we will make the
two following assumptions on our random roughness:\\
{\em (H1)  The measure $\pi$ is supported by 
$$ P_\alpha = \left\{ \omega : \R \mapsto (-1,0), \: \| \omega
  \|_{C^{2,\alpha}} \: \le  K_\alpha \right\} $$
for some $\alpha > 0$ and some $K_\alpha > 0$.} \\
{\em (H2) The randon boundary has no correlation at large distances, that
  is the $\sigma$-fields
$$\sigma\left( s \mapsto \omega(s), \: s \le a \right)  \: \mbox{ and } \: 
\sigma\left( s \mapsto \omega(s), \: s \ge b \right)$$ 
are independent for $b-a \ge \kappa$, for some $\kappa > 0$.}
 \\
Under these assumptions, the main theorem of the paper reads:
\begin{theorem} \label{refinedestimate}
For small enough $\phi$ and under (H1)-(H2),
 the following refined estimate holds:
$$   \| u^\eps - u^N \|_{L^2_{uloc}(P \times \Omega)} = 
O(\eps^{3/2} |\ln(\eps)|^{1/2}). $$
\end{theorem}
Before entering the proof of this theorem, let us give a few hints.
 Theorem \ref{refinedestimate} is deduced from  a central limit theorem for
 the quantity $v(\omega,y)-(\alpha,0)$. Broadly, this theorem comes from
 good properties of the random variables
$$ X^n(\omega) \: = \: \int_n^{n+1} v(\omega,y_1, 0) \, dy_1. $$
Due to the elliptic nature of the Stokes operator, such random variables
are not independent. However, under assumption (H2), we are able to prove
that the correlation terms $E(X_n \, X_0)$ decay fast enough as $n
\rightarrow \infty$. As a result,
one can prove a central limit theorem on $X_n$, and then a similar one on
$v - (\alpha,0)$. We point out that such type of results  for dependent
variables with strong  decay of correlations is quite classical and has
been used in various fields. We refer to \cite{Baladi:2001} for a review
paper related to dynamical systems, and to recent articles \cite{Varadhan:2007,
  Bouard:2007}  for applications  in a  PDE context. 

As a consequence of this central limit theorem,
 we show that the boundary layer converges to a constant as
 $|y_2^{-1/2}|$. Note that this is in sharp contrast with the periodic
 case, where exponential convergence holds (we stress that  periodic
 boundaries are highly correlated, thus far from satisfying (H2)).This
 speed of convergence is resposible for the $\eps^{3/2} |\ln(\eps)|^{1/2}$
 in the  Navier wall law.  

 The main difficulty  is to obtain the decay of correlations
 of variables like $X_n$. The proof relies on precise estimates of the Green
 function for the Stokes operator above a non flat  boundary. Such
 estimates follow from sharp elliptic regularity results, 
where one must pay attention to the oscillation of the 
boundary.  This is achieved under the regularity assumption (H1), using
 ideas of Avellaneda and Lin for homogenization of elliptic systems
 \cite{Avellaneda:1987, Avellaneda:1991}.

\section{Boundary layer decay and Navier approximation}
In this section, we explain how Theorem \ref{Navestimates} follows
from estimates on the solution $v$ of \eqref{BL}. Such estimates will be
established in the following sections. At first, we remind the main
features of $v$, as stated in article \cite{Basson:2007}. 

\subsection{The boundary layer system}
As emphasized in the introduction,  to solve 
\eqref{BL} in a deterministic way, that is for each possible boundary
$\omega$, is still unclear. Hence, one must take advantage of the
probabilistic setting. 
First, notice that a reasonable solution $v$ should satisfy: 
\begin{equation} \label{stationarity}
 v(\tau_h(\omega), y_1, y_2) \: = \: v(\omega, y_1+h, y_2). 
\end{equation}
Together with the stationarity assumption, this relation sort of
substitutes to the identity 
$$ v(y_1 + L, y_2) \: =  \: v(y_1, y_2) $$
used in the treatment of $L-$periodic roughness. It allows to extend the
well-posedness result, through an appropriate variational formulation. 

\medskip
 This formulation  has been described in article \cite{Basson:2007}.
 First, one introduces the new unknown 
$$ w(y) \: :=  \: v(y) \: + \: (y_2,0) \, \mathbf{1}_{\{y_2 < 0\}}(y), $$ 
and replace system \eqref{BL} by 
\begin{equation} \label{BL2}
\left\{
\begin{aligned}
& -\Delta w + \na q = 0, \: x \in \Omega^{bl}\setminus \{ y_2 = 0\}, \\
& \div w = 0, \: x \in \Omega^{bl},\\
& w\vert_{\pa \Omega^{bl}} \: = \: 0,\\
& [ w ]\vert_{y_2 = 0} = 0, \quad [\pa_2 w - (0,q)]\vert_{y_2 = 0} =  (-1,0),  
\end{aligned}
\right.
\end{equation}
where $[ \cdot ]\vert_{y_2 = 0}$ denotes the jump at $y_2=0$. Then, one
multiplies formally  the Stokes equation by a  test
function $w' = w'(\omega,y)$ that satisfies $\div w' = 0$, $w'\vert_{\pa
 \Omega^{bl}} = 0$. Integrating by parts over $\Omega^{bl} \cap \{ |y_1 <
1\}$  yields 
 $$ \int_{\Omega^{bl} \cap \{ |y_1 < 1\}} \na w \cdot \na w' \: = \:
   \int_{\{ |y_1| < 1, \: y_2 = 0\}}  w'_1 \: + \int_{\Omega^{bl} \cap
 \{|y_1| = 1\}}  \left( \pa_n w - q n \right) w'. $$
Finally, if $w$, $w'$ satisfy relation \eqref{stationarity}, one can
 integrate with respect to $\omega$,  and thanks to the  stationarity of 
$\pi$, get rid of the annoying boundary term at the r.h.s: 
$$ \E  \int_{\Omega^{bl} \cap \{ |y_1 < 1\}} \na w \cdot \na w' \: = \:
   \E \int_{\{ |y_1 < 1, \: y_2 = 0\}}  w'_1 $$

Afterwards, this formal variational formulation can be rigorously defined
 and solved: in short, one can apply the Riesz theorem in a functional
 space of Sobolev type, made of functions $w$ such that  
$$ \E \int_{\Omega^{bl} \cap \{  |y_1 < 1\}} | \na w |^2 \: < \: +\infty, $$
and satisfying almost surely \eqref{stationarity}, together with
$ \div w = 0, \: w\vert_{\Omega^{bl}} = 0$. 
We refer to \cite{Basson:2007}
 for all details. Note that stationarity implies: 
$$ sup_{t, R} \: \E \: \frac{1}{R} \int_{\Omega^{bl} \cap 
\{  |y_1 -t | < R\}} |
 \na w |^2 \: = \:  \E \int_{\Omega^{bl} \cap \{  |y_1 < 1\}} | \na w |^2 
\: < \: +\infty. $$

Back to the original system \eqref{BL}, this variational solution $w$
 provides almost surely a solution $v(\omega, \cdot) \in
 H^1_{loc}\left(\Omega^{bl}\right)$ in the sense of
 distributions. Moreover, the ergodic theorem yields (see
 \cite{Basson:2007})   
$$ \sup_{R} \frac{1}{R} \int_{\Omega^{bl} \cap \{  |y_1| < R\}}  |  \na v
 |^2 \: < \:  +\infty, \: \mbox{ almost surely}. $$

\medskip
In order to understand the origin of the Navier approximation, the next
 step is to
 describe the behavior of $v$ as $y_2 \rightarrow +\infty$. 
 For periodic roughness, one can show
exponential convergence of $v$ to a constant vector field $v^\infty =
(\alpha, 0)$. However, the rate of convergence goes to zero with the period
$L$. When dealing with stationary random boundaries, that broadly speaking
contain all periods, the exponential decay does not hold {\it a priori}. In
other words, there is a problem associated to the Fourier spectrum, that is
discrete in the periodic case, and may accumulate to zero in the random
case. 

Again, this problem has been (partially) overcome in \cite{Basson:2007}.
 The first
step is to obtain a  representation of  $v$ in terms of a Stokes double
layer potential, {\it cf} proposition ??. Almost  surely, for any  
\begin{equation*} 
\begin{aligned}
 v(\omega, y) \: & = \:  \int_\R G(t,y_2) \, v(\omega,y_1-t,0) \, dt \\
                 & = \: -\int_\R t \, \pa_t G(y_2) \frac{1}{t} \int_0^t
                 v(\omega,y_1-s,0) \, ds \, dt 
\end{aligned}
\end{equation*}
where $G$ is the Poisson type kernel for the Stokes operator over a half
space. Then, the ergodic theorem and a few calculations yield:   
$$   \frac{1}{t} \int_0^t v(\omega,y_1-s,0) \, ds \,  \: \rightarrow
v^\infty(\omega) \: = \: (\alpha(\omega), 0), \quad t \rightarrow \pm \infty,$$
where the convergence holds almost surely (locally uniformly in $y_1$), as
well as in $L^2(P)$ (uniformly in $y_1$). In the case where the
stationary measure $\pi$ is ergodic, the constant $\alpha$ does not depend
on $\omega$.
Finally, back to the integral representation, and with 
similar treatment for  derivatives of $v$: 
\begin{equation} \label{softestimbl}
  \forall \beta \in \N^2, \: |\beta| \ge 1, \quad  
 \E \left|v(\cdot, 0,,y_2) - \alpha(\cdot)\right|^2  \:  + \: 
y_2^{2|\beta|}   \E
  \left|\pa^\beta_y \,v(\cdot, 0,y_2)\right|^2 \: \xrightarrow[y_2 \rightarrow
    +\infty]{} \: 0.
\end{equation}
We refer to  \cite[Proposition 13]{Basson:2007} for all details. 

\subsection{Refined estimate for Navier wall law}
Most of the analysis of the present paper  will be devoted to a  refined 
 asymptotic estimate of the boundary layer:
\begin{theorem} \label{theosharp}
Under assumptions (H1), (H2), for all $\beta \in \N^2$,  
\begin{equation} \label{sharpestimbl}
 y_2^{2 |\beta|+1} \E \left|\pa^\beta \left( 
v(\cdot, 0,y_2) - \alpha(\cdot)\right) \right|^2  \: \xrightarrow[ y_2
  \rightarrow +\infty]{} \sigma_\beta \ge 0. 
\end{equation}
\end{theorem}
It is of course a much sharper convergence result than \eqref{softestimbl}. 
Before tackling its proof, we explain how it implies theorem
\ref{Navestimates}. Arguments are direct  adaptation from section 5  in
\cite{Basson:2007}.  

\medskip
On the basis of the boundary layer analysis, one can  build an
approximation of $u^\eps$ of boundary layer type. Namely, we introduce
$$ u^\eps_{app}(\omega,x) \: = \:  u^0(x) \: + \: 
\, 6 \, \phi\,\eps \,  v\left(\omega,\frac{x}{\eps}\right) +
 \, 6 \, \phi\,\eps \,   u^1(\omega,x) +
\, 6 \, \phi\,\eps \,r^\eps(\omega,x). $$
In this approximation, $u^0$ is the Poiseuille flow and $v(\omega,\cdot)$
 is the boundary layer solution of \eqref{BL}. As $v$ does not
converge to zero at infinity, we add a large scale corrector $u^1$
satisfying: 
\begin{equation} \label{Couette}
\left\{
\begin{aligned} 
& u^0 \cdot \na u^1 + u^1 \cdot \na u^0 -\Delta u^1 + \na p = 0, \: x \in
  \Omega,\\ 
& \div u^1 \: = \:  0, \: x \in \Omega,\\
& \int_0^1 u^1 \cdot e_1  dx_2  = 0, \\ 
& u^1\vert_{y_2 = 0} = 0, \quad u^1\vert_{y_2 = 1} = -(\alpha, 0). 
\end{aligned}
\right.
\end{equation}
It is just a combination of a Couette and a Poiseuille flow:
$u^1  =   \alpha  x_2 \left( 2 - 3 x_2 \right) e_1$.
 Still, this approximation does not vanish at the boundary, which explains
 the addition of another term  $r^\eps(\omega,x)$. It must satisfy 
\begin{equation} 
\left\{
\begin{aligned}
& r^\eps(x_1,0) \: = \: 0, \\
& r^\eps(x_1,1) \: = \: v\left(\frac{x_1}{\eps}, \frac{1}{\eps}\right) -
  (\alpha,0),  \\
& \div r^\eps \: = \:  0, \: x \in \Omega.   
\end{aligned}
\right.
\end{equation}
This remainder can be taken  small in the sense of 
\begin{proposition} 
This problem possesses a (non unique) solution $r^\eps$ such that
$$ \sup_x \, \E  \: \| r^\eps \|^2_{H^2(B(x,1) \cap \Omega)} \: = \:
O(\eps | \ln \eps|). $$ 
\end{proposition}
{\em Proof:} The proof of this result mimics the one of proposition 14 in
\cite{Basson:2007}.
 The corrector $r^\eps$ can be chosen  in the  form 
$$ r^\eps = \na^\bot \psi, \quad \psi = a(x_1)x_2^3  + b(x_1) x_2^2 +
 c(x_1) x_2 + d(x_1).$$
The streamfunction $\psi$ is determined up to a constant and 
 polynomial in $x_2$. Its coefficients  have explicit dependence on
 $v-(\alpha,0)$.  Hence, the $H^2$ estimate on $r^\eps$ follows from the
control of various terms involving $v-(\alpha,0)$.  For instance, one must
bound the $L^2(P \times (-1,1))$ norm of  
\begin{equation*}
\begin{aligned}
 \int_0^{x_1} v_2\left(\omega, \frac{t}{\eps},\frac{1}{\eps}\right) dt
 \: & = \: \eps \int_0^{x_1/\eps} 
v_2\left(\omega,y_1,\frac{1}{\eps}\right) dy_1 \\
&  = \:  \eps \int_{\omega(x_1/\eps)}^{1/\eps} 
( v_1 - \alpha)\left(\omega,\frac{x_1}{\eps}, y_2 \right) \,  dy_2 \\
& \: - \:  \eps 
\int_{\omega(0)}^{1/\eps} (v_1 - \alpha)(\omega,0,y_2) \, dy_2  \: :=
I^\eps(\omega,x_1) \: - \:  I^\eps(\omega,0)  
\end{aligned}
\end{equation*}
where the last equality comes from the Stokes formula. Using stationarity
of $\pi$, we get 
$$ \E \left\|x_1 \mapsto  \int_0^{x_1} v_2\left(\cdot, \frac{t}{\eps},
\frac{1}{\eps}\right) dt  \right\|_{L^2(-1,1)}^2
 \: \le \: 4 \E | I^\eps(\cdot,0) |^2. $$
Thanks to the refined estimate \eqref{sharpestimbl}, we finally obtain 
\begin{align*}
\E | I^\eps(\cdot,0) |^2 \: & \le \:  C \left( 
\eps^2 \,  \E \int_{\omega(0)}^{1}  |(v_1 -
\alpha)(\cdot,0,y_2)|^2  dy_2 \: + \:  \eps \int_1^{1/\eps} \E \, |(v_1 -
\alpha)(\cdot,0,y_2)|^2  dy_2 \right)\\
&  \le \:  C' \eps^2 + C'' \eps \int_1^{1/\eps} y_2^{-1} \, dy_2 = O(\eps
|\ln \eps|)
\end{align*}
All other terms involve similar computations. These are straightforwardly
adapted from the proof of proposition 14 in \cite{Basson:2007}, using
\eqref{sharpestimbl} instead of \eqref{softestimbl}. 

\medskip  
Once the approximate solution $u^\eps_{app}$ is built, one can obtain by
energy estimates the following bounds, for $\phi$ small enough:
\begin{align*}
& \| u^\eps - u^\eps_{app} \|_{L^2_{uloc}(P \times \Omega)} \: = \:
O\left(\eps^{3/2} |\ln(\eps)|^{1/2}\right), \\
& \| u^\eps_{app} - u^N \|_{L^2_{uloc}(P \times \Omega)} \: = \:
O\left(\eps^{3/2} |\ln(\eps)|^{1/2}\right),
\end{align*}
which of course imply theorem \ref{Navestimates}. As the proof is very similar
to what was done in paper \cite{Basson:2007},
 we do not expand more and refer to it for all details. 

\section{A central limit theorem}
Up to the end of the paper, we will assume (H1)-(H2), and focus on
 theorem \ref{theosharp}. It is classical that (H2) implies ergodicity of
 $\pi$, so that the constant $\alpha$ does not depend on $\omega$. 
 We start again from an integral representation
\begin{equation} \label{layer}
\begin{aligned}
 \pa^\beta_y \left( v(\omega, 0,y_2) - (\alpha, 0) \right)  \: & = \:  
\int_\R \pa^{\beta_1}_t \pa^{\beta_2}_{y_2} \, 
 G(t,y_2) \, \left(v(\omega,-t,0) - (\alpha,0) \right)  \, dt \\
& = \: -\int_\R  \, \pa^{\beta_1+1}_t \pa^{\beta_2}_{y_2} \, 
G(t,y_2) \int_0^t
\left(v(\omega,-s,0) - (\alpha,0) \right) ) \, ds \, dt,
\end{aligned}
\end{equation}
where the matrix kernel $G$ is given by 
$$ G(y) \: = \: \frac{2 y_2}{\pi (y_1^2 + y_2^2)^2} \begin{pmatrix} y_1^2 &
  y_1 y_2 \\ y_1 y_2 & y_2^2 \end{pmatrix} $$
We introduce 
$$V(\omega,t) \: := \:
 \int_{0}^t \left(v(\omega,-s,0) - (\alpha,0) \right)  \, ds.  $$
A simple change of variable leads to 
\begin{equation}
  y_2^{|\beta|+1/2} 
\pa^\beta_y \left( v(\omega, 0,y_2) - (\alpha, 0) \right)  \:
 = \:   \int_\R \pa^{\beta_1+1}_t \pa^{\beta_2}_{y_2} G(u,1) \,
y_2^{-1/2} \,  V(\omega, y_2 u) \, du. 
\end{equation}
Our first goal is to show that the l.h.s. converges in law  
to a gaussian distribution, for all $\beta$. We will focus on the case
  $|\beta| = 0$, the other cases being handled in the exact same way.
We state 
\begin{proposition} \label{tcl}
The function $V$ satisfies the following properties:
\begin{description}
\item[i)] $\E \left|V(\cdot,t)\right|^2 \: \le \: C  |t|$ 
\item[ii)]   The random process   $y_2^{-1/2} \, V(\omega, y_2 \, u)$
 converges weakly to a gaussian process $B(\omega,u)$  as $y_2$ goes to 
 infinity.
\item[iii)]   The covariance matrices also converge, that is for all
  indices $i,j$  and for all $s,t$, 
$$\E \, y_2^{-1} \, V_i(\cdot, y_2 \, s) \, V_j(\cdot, y_2 \, t)
  \xrightarrow[y_2 \rightarrow +\infty]{} \E B_i(\cdot, s)
   B_j(\cdot, t)$$
\end{description}
\end{proposition}
We remind that the process $X^n(\omega,t)$ with values in $\R^2$
 {\em converges weakly} to
$X(\omega,t)$ if, for all $T > 0$ and all continuous bounded
function ${\cal F} : C\left( [-T,T], \R^2 \right) \mapsto \R$, 
$$   \E  {\cal F}(X^n)  
\: \xrightarrow[n \rightarrow +\infty]{} \: \E {\cal F}(X). $$
Theorem \ref{theosharp} is then a direct consequence of
\begin{corollary}
The random process  $y_2^{1/2} \left(v(\omega,0,y_2) -
(\alpha,0)\right)$  converges in law  to a gaussian vector with zero average. 
Moreover, for all $i,j$, 
$$ \E  \left(v_i(\cdot,0,y_2) - (\alpha,0)_i\right) \,
\left(v_j(\cdot,0,y_2) - (\alpha,0)_j\right) \xrightarrow[y_2 \rightarrow
  +\infty]{} \sigma_{ij}, $$
where $\sigma$ is the covariance matrix of this  gaussian vector. 
\end{corollary}
{\em Proof of the corollary:} To prove convergence in law to a gaussian
 vector  ${\cal N}_{\sigma}$ of covariance matrix $\sigma$, we need
 to show that for any $F \in C^\infty_c(\R)$, 
$$\E \, F \left( \int_\R \pa_t G(t,1) \, 
y_2^{-1/2} \, V(\cdot, y_2 t) \, dt\right)   \xrightarrow[y_2 \rightarrow
  +\infty]{} \E \, F\left({\cal N}_{\sigma}\right). $$
 Unsurprisingly, we take
$$  {\cal N}_{\sigma} \: := \: \int_\R \pa_t G(t,1) B(\omega,t) \, dt. $$ 
We decompose, for any  $T > 0$, 
\begin{align*}
& \E \, F \left( \int_\R \pa_t G(t,1) \, 
y_2^{-1/2} \, V(\cdot, y_2 t) \, dt\right) - \E \, F\left({\cal
  N}_{\sigma}\right)  \\
\: =  &\: \E \, F \left( \int_{-T}^T \pa_t G(t,1) \, 
y_2^{-1/2} \, V(\cdot, y_2 t) \, dt\right) - 
\E \, F \left( \int_{-T}^T \pa_t G(t,1) \,  B(\cdot,t) \, dt\right)  \\
 \: + & \:\E \, F \left( \int_\R \pa_t G(t,1) \, 
y_2^{-1/2} \, V(\cdot, y_2 t) \, dt\right) -   
\E \, F \left( \int_{-T}^T \pa_t G(t,1) \, 
y_2^{-1/2} \, V(\cdot, y_2 t) \, dt\right)  \\
\: + \:& \E \, F \left( \int_\R \pa_t G(t,1) \, 
B(\cdot,t) \, dt\right) -   
\E \, F \left( \int_{-T}^T \pa_t G(t,1) \, 
 B(\cdot,t) \, dt\right) \\
 \: :=  &\: J_1 + J_2 + J_3 
\end{align*}
We now show that expressions $J_1$, $J_2$,   converge to zero
($J_3$ is similar to $J_2$ and simpler). Let
$\delta > 0$. We have
\begin{align*}
 |J_2| \: & \le \: \max |F'|   \: \E \int_{|t|>T}
 \left|\pa_t G(t,1)\right| \, 
\left|y_2^{-1/2} \, V(\omega, y_2 t) \right| \, dt \\
& \le \: \max |F'|  \: \left( \int_{|t|>T} \left|\pa_t G(t,1)\right| dt
\right)^{1/2} \, \left( \int_{|t|>T}  \left|\pa_t G(t,1)\right| 
\, \E  \left(  y_2^{-1} \,  |V(\omega, y_2 t)|^2 \,\right) 
 dt \right)^{1/2} \\
& \le \: C \left( \int_{|t|>T} \left|\pa_t G(t,1)\right| dt
\right)^{1/2} \left( \int_{|t|>T} \left|\pa_t G(t,1)\right| 
t \, dt \right)^{1/2} 
\end{align*}
where the last line comes from  point ii) of  proposition \ref{tcl}. Thus,
for $T$ large enough, independently of $y_2$, $|J_2| \le \delta/2$. Such
$T$ being fixed, for $y_2$ large enough, we get  
$|J_1| \le \delta/2$ by point i) of proposition \ref{tcl}. This yields 
convergence in law. The convergence of the covariance matrix 
\begin{align*}
 & \E  \left(v_i(\cdot,0,y_2) - (\alpha,0)_i\right) \,
\left(v_j(\cdot,0,y_2) - (\alpha,0)_j\right) \\
& = \:
 \int_\R \int_\R \E \left( \pa_t G(s,1) 
y_2^{-1/2} \, V(\cdot, y_2 \, s) \right)_i  \, \left( \pa_t G(t,1) 
y_2^{-1/2} \, V(\cdot, y_2 \, t) \right)_j  \, ds \, dt  
\end{align*}
follows from the dominated convergence theorem, using i) and iii) of
proposition \ref{tcl}. We get 
$$ \sigma_{ij} \: = \: \E \int_\R \int_\R \left( \pa_t G(s,1) B(\cdot,s)
\right)_i \, \left( \pa_t
G(t,1) B(\cdot,t) \right)_j \, ds \, dt. $$
This concludes the proof of the corollary.  

\medskip
It remains to prove theorem \ref{tcl}. Note that point ii) 
  is essentially  a central limit theorem for the sequence of random
  variables 
$$X^n(\omega) = F \circ \tau_n (\omega), \quad F(\omega) =  \int_0^1
  \left(v(\omega,t,0) - (\alpha(\omega),0) \right) \, dt.$$ 
The problem is  that these random variables are not
  independent, due to  ``propagation of information at infinite speed''
 in the Stokes  system. To establish a central limit theorem for such type of
  sequences is a classical  question. The basic idea is that
  one can extend the central limit theorem to non independent
  sequences that feature  a good decay of correlations
  as $n$ goes to infinity. We now illustrate this general principle
  on our problem, using assumption (H2). We follow the presentation of  
article \cite{Varadhan:2007}, in which a
  similar question arises for a  semilinear heat equation
  with random source. Let ${\cal C}^n$ the $\sigma$-algebra generated
  by the applications $y_1 \mapsto \omega(y_1)$, $|y_1| < n$.We state the
  following lemma: 
\begin{lemma} \label{lemmeUn}
Suppose that $\: v^n \: := \: \E\left( v(\cdot,0,0)  \: |  \:  {\cal C}^n\right)$  satisfies 
$$\E \left| v^n - v(\cdot,0,0) \right|^2  \: \le \: 
C \, n^{-\alpha}$$
 for some $\alpha >1$. Then, proposition \ref{tcl} holds. 
\end{lemma}  
{\em Proof of the lemma:} We  write the decomposition
$$  v(\cdot,0,0) - (\alpha,0) \: = \: 
v^1 - (\alpha,0) \: + \: \sum_{j=1}^{+\infty} \left( v^{2^j}
- v^{2^{j-1}} \right) $$ 
with the sum converging in $L^2(P)$. The corresponding sum  for $V$ is 
$$ V \: = \: \sum_{j=0}^{+\infty} \,  V^{j}, \quad V^j(\omega,t) =
\int_{0}^{t} \left( v^{2^j}
- v^{2^{j-1}} \right)\circ \tau_s(\omega)  \, ds, $$
where $v^{1/2} := (\alpha,0)$.   
Then, we have:  $ \: \| V(\cdot,t) \|_{L^2(P)} \, \le \,
\sum_{j=0}^{+\infty} \: \|  V^{j}(\cdot,t) \|_{L^2(P)}$. By the
assumption of independence at large distances, 
the correlations $\E \left( v^{2^j} \circ 
\tau_t(\omega) \, v^{2^j} \circ \tau_{t+s}(\omega) \right)$ and 
$\E \left( v^{2^j} \circ 
\tau_t(\omega) \, v^{2^{j-1}} \circ \tau_{t+s}(\omega) \right)$
vanish
for $|s| \ge \kappa + 2^{j+1}$. We introduce
 $$n \: :=  \:  \left[ |t|/(\kappa+2^{j+1})\right].$$
If $n=0$, we just write
$$   \E  \left|V^{j}(\cdot,t)\right|^2  \: \le \: |t|^2 \E
  \left|v^{2^j} - v^{2^{j-1}} \right|^2. $$
If $n \ge 1$, we decompose 
\begin{equation*}
\begin{aligned}
 \E  \left|V^{j}(\cdot,t)\right|^2 \:  & =  \: \E \left| \sum_{k=0}^{n-1}
  \int_{kt/n}^{(k+1)t/n} \left( v^{2^j} \circ \tau_s - v^{2^{j-1}}
  \circ \tau_s \right)  \, ds  \right|^2  \\
&   \le \:  \E \left| \int_0^{t/n}  \sum_{k=0}^{n-1} \left( v^{2^j} \circ
  \tau_{s+kt/n} -  v^{2^{j-1}} \circ  \tau_{s+kt/n} \right)  \, ds
  \right|^2  \\
&  \le  \:  2 \left(\kappa + 2^{j+1} \right)
 \, \int_0^{|t|/n}  \sum_{k=0}^{n-1} \E \left| v^{2^j} -
  v^{2^{j-1}} \right|^2 
\end{aligned}
\end{equation*} 
  Using the bound on the
conditional expectations, we end up with
\begin{equation} \label{estimateVj}
\| V^{j}(\cdot,t) \|_{L^2(P)}^2  \: \le \: C \, |t| \: \min(|t|, 2^j) \, 
2^{-j \alpha}.
\end{equation}  
for some constant $C = C(\kappa)$. We thus get i). To prove ii), we just write
the decomposition
$$  v(\cdot,0,0) - (\alpha, 0) \: = \: 
 \sum_{j=0}^{+\infty} \left( v^{2^j} \: - \: v^{2^{j-1}} \right), \quad 
 y_2^{-1/2} \,V(\omega, t y_2) \: = \: y_2^{-1/2} \sum_{j=0}^{+\infty}
 \,V^j(\omega, t y_2).  $$  
It is well-known that each finite sum satisfies a central limit theorem,
 that is  
$$ \forall k, \quad  y_2^{-1/2} \sum_{j=0}^{k}  \,V^j(\omega, t y_2)
 \xrightarrow[y_2 \rightarrow +\infty]{} B^k(\omega,t) $$
in the sense of weak convergence, to some gaussian process $B^k(\omega,t)$.
 Moreover, the covariance matrix also converges, that is 
 $$ y_2^{-1} \E   \sum_{j=0}^{k}  \,V^j_l(\cdot, s y_2)  \,
 \sum_{j=0}^{k}  \,V^j_m(\cdot, t y_2)  \: \xrightarrow[y_2 \rightarrow
 +\infty]{} 
 \E  B^k_l(\omega,s) \, B^k_m(\omega,t). $$
In short, it is due to the fact that the random variables 
$$ X^{n,j}(\omega) =   F^j \circ \tau_n(\omega), \quad
  F^j(\omega) =   \int_0^1 \left( v^{2^j} -  v^{2^{j-1}} \right) \circ
  \tau_t(\omega) \, dt, \quad 
  n \in \Z,$$ 
have zero correlations at large distances: see \cite[theorem (7.11)
  p424]{Durrett:1996}  for
  a similar result and detailed proof.
 Moreover, thanks to estimate \eqref{estimateVj}, the  remainder 
$$ R^k(\omega,t,y_2) \: = \:  \sum_{j=k}^{+\infty} \, y_2^{-1/2}
 \,V^j(\omega, t y_2) $$ 
converges to zero as $k \rightarrow +\infty$, locally uniformly in
$t$, uniformly in $y_2$. Hence, points (ii) and (iii) of proposition \ref{tcl}
hold, which ends the proof of the lemma.  

\medskip
We still have to estimate the variance of $v^n - v(\cdot,0,0)$. 
Following \cite{Varadhan:2007}, we can turn this question into a question
of domain of dependence for solutions of \eqref{BL}. Precisely,
starting from the measure $\pi$ on $P$, we  define a new measure
$\pi^n$ on the product space 
$$P^n \: = \: \left\{ (\omega_1, \omega_2) \in P \times P, \quad
 \omega_1(t) = \omega_2(t), \: |t| \le n \right\}. $$  
endowed with its cylindrical $\sigma-$field. Namely, $\pi^n$ is  defined in
 the following way:  
\begin{enumerate}
\item $\pi^n(A \times A) \, :=  \, \pi(A), \: \forall A \in {\cal C}^n$, which 
  determines $\pi^n$ over  the sub $\sigma-$field ${\cal D}^n$ generated by the
  applications $\displaystyle t  \mapsto (\omega_1(t), \omega_2(t)), \: |t| \le n$.  
\item For all $k \ge 1$, for all $t^1, \dots, t^k$ with $|t^j| > n$, for
  all  borelian subsets  $B^1_1, \dots, B^k_1$, $B^1_2, \dots, B^k_2$  of
  $\R$, and  
$$A_1 \: := \: \cap_{j=1}^k \{ \omega_1, \:  \omega_1(t_j) \in
 B_1^j \}, \quad A_2 \: := \: \cap_{j=1}^k \{ \omega_2, \:  \omega_2(t_j) \in
 B_2^j \}$$
$\pi^n(A_1 \times A_2 \, | \,  {\cal D}^n)(\omega_1,\omega_2) \, :=
 \,\pi(A_1 \, | \, {\cal 
 C}^n )(\omega_1) \: \pi(A_2 \, | \, {\cal C}^n )(\omega_2)$, which
 determines $\pi^n$ 
 conditionaly to ${\cal D}^n$.    
\end{enumerate}
It is then easy to derive the following identity, see \cite{Varadhan:2007}:
$$ \E \left|v^n - v(\cdot,0,0 )  \right|^2 \: = \: \frac{1}{2} \int_{P^n} \left| 
v(\omega_1,0,0) - v(\omega_2,0,0) \right|^2 d\pi_n. $$
Thus, if  $\Omega^{bl}(\omega_1)$ and $\Omega^{bl}(\omega_2)$ are two
boundary layer domains with boundaries that coincide over $[-n,n]$,  we need to estimate the
difference of the corresponding boundary layer solutions
$v(\omega_1,0,0)$ and  $v(\omega_2,0,0)$. This is the purpose of the
next section.

\section{Decay of correlations}
{\em Throughout the rest of the paper, we will assume (H1).} 
The main difficuly is to prove the following 
\begin{proposition} \label{propdecay}
Under assumption (H1), for all $0 < \tau < 1$, for almost every  
$\omega_1, \omega_2 \in P$,  
\begin{equation} \label{decaybound}
\left| v(\omega_1,0,0) - v(\omega_2,0,0) \right| \: \le \: 
\frac{C}{n^{2\tau-1}},
\end{equation}
where $C$ does not depend on $\omega_1, \omega_2$.
\end{proposition} 
Together with the results of the preceding section, this proposition  
concludes the proof of
theorem \ref{theosharp} (take $\tau > 3/4$),
 and therefore the proof of the main theorem \ref{Navestimates}. In fact, the
sharper bound  
$$\left| v(\omega_1,0,0) - v(\omega_2,0,0) \right| \: \le \: \frac{C}{n},$$
that is with $\tau=1$ would still be true. We will discuss this
 briefly in the last
section of the paper. For the sake of brevity, we  only prove here the
weaker form \eqref{decaybound}.
    
The main difficulty  is that the boundary layer
solutions $v(\omega_1,y)$ and $v(\omega_2,y)$ of \eqref{BL} are not
defined on the same domain, so that  estimates on the difference
are not directly available. If the Poisson equation rather than the  Stokes
system  was considered,  representation of the solution in
terms of Brownian motion would allow to conclude quite easily. 
Again, this  will be explained in the last section of the paper. 
 
\medskip
In the case of system \eqref{BL}, we are not aware of such representation,
and the bound \eqref{decaybound} will come from an accurate description of
the (matrix) Green function of the Stokes operator above a humped boundary. 
We consider for all $\omega \in C^{2,\alpha}$, and for all $z \in 
\Omega^{bl}(\omega ) \: = \: \{  y_2 > \omega(y_1)\}$, the system:
\begin{equation} \label{green}
\left\{
\begin{aligned}
 -\Delta G_\omega(z,\cdot) + \na P_\omega(z,\cdot) = \delta_z \, 
I_2 & \quad \mbox{ in
} \Omega^{bl}(\omega), \\
 \div G_\omega(z,\cdot) = 0 & \quad \mbox{ in } \Omega^{bl}(\omega), \\
G_\omega(z, \cdot) = 0 & \quad \mbox{ on } \pa \Omega^{bl}(\omega). 
\end{aligned}
\right.
\end{equation}
where $\delta_z$ is the Dirac mass at $z$, and $I_2$ is the $2\times2$
identity matrix. Let us  remind how to build  the
matrix Green function $(G_\omega, P_\omega)$. Up to a vertical translation
of the 
domain, we can first assume that $z_2 > 0$. We then introduce the Green
function $(G_0, P_0)$ for the Stokes operator in the
upper-half plane, see \cite{Galdi:1994}.  Extending $\displaystyle 
G_0(z, \cdot), P_0(z, \cdot)$ by  $0$ for $y_2 < 0$, the
functions 
$$\displaystyle H(z, \cdot) \: :=  \; G_\omega(z, \cdot) - G_0(z, \cdot),
\quad Q(z, \cdot) 
\: := \: P_\omega(z, \cdot) - P_0(z, \cdot)$$ 
 satisfy formally 
\begin{equation*} 
\left\{
\begin{aligned}
 -\Delta H(z,\cdot) + \na Q(z,\cdot) = 0 & \quad \mbox{ in
} \Omega^{bl}(\omega), \\
 \div H(z,\cdot) = 0 & \quad \mbox{ in } \Omega^{bl}(\omega), \\
H_\omega(z, \cdot) = 0 & \quad \mbox{ on } \pa \Omega^{bl}(\omega), \\
 \left[ H(z,\cdot)\right] = 0, & \quad \bigl[ \pa_2 H(z,\cdot) - Q(z, \cdot)
   \otimes  \,e_2 \bigr] = - \left[ \pa_2 G_0(z,\cdot) - P_0(z, \cdot)\,
    \otimes e_2  \right], 
\end{aligned}
\right.
\end{equation*}
where $\left[ \: \cdot \:  \right]$ is the jump along  $\{ y_2 = 0 \}
\cap \Omega^{bl}(\omega)$. The jump on the derivative is explicit, as 
$$ \pa_2 G_0(z,(y_1,0^+)) - P_0(z,(y_1,0^+)) \otimes e_2 \: = \: 
\frac{2 z_2}{\pi ((z_1-y_1)^2 + z_2^2)^2} \begin{pmatrix} (z_1-y_1)^2 &
  (z_1-y_1)  y_2 \\ (z_1-y_1)  y_2 & y_2^2 \end{pmatrix}.
$$ 
Standard variational formulation  yields existence and uniqueness of a solution
$H(z,\cdot)$ with $\na H(z, \cdot)$ in  $L^2$. In turn, this provides a
unique solution $G_\omega(z, \cdot)$ to \eqref{green}. The corresponding
pressure field $P_\omega(z, \cdot)$ is determined up to the addition of a
constant matrix. Note that uniqueness yields the relation
\begin{equation} \label{stationG}
 G_{\tau_h(\omega)}(z,y) \: = \: G_\omega((z_1+h,z_2), (y_1+h,y_2)). 
\end{equation}

\medskip
Our  key estimate is provided by 
\begin{lemma}  \label{green2}
For all $0 < \tau < 1$, for all $z,y \in \Omega^{bl}(\omega)$ satisfying
$|z-y| \ge 1$,  we have 
\begin{equation} \label{greenestim1}
\sum_{|\beta| \le 2} |\pa^\beta_y G_\omega(z,y)| \: + \: |\na_y
 P_\omega(z,y)| \:  \le 
 \: C \,  \frac{ \delta(z)^\tau \, (1+ \delta(y))^\tau}{|z-y|^{2\tau}}, 
\end{equation}
where $\delta(\cdot)$ denotes the distance to the boundary $\pa
 \Omega^{bl}(\omega)$, and $C$ is  a constant
  depending only on  $\tau$ and on $\| \omega \|_{C^{2,\alpha}}$.
\end{lemma}
Note that by symmetry of $G$, we also have 
\begin{equation} \label{greenestim2}
\sum_{|\beta| \le 2} |\pa^\beta_z G_\omega(z,y)| \: \le \: C \,  \frac{
  (1+\delta(z))^\tau \, (1+ \delta(y))^\tau}{|z-y|^{2\tau}}. 
\end{equation}
Moreover, in the course of the proof of lemma \ref{green2}, we will show
that for all $y,z \in \Omega^{bl}(\omega)$, 
\begin{equation} \label{greenestim3}
 | G_\omega(z,y) | \: \le \: C \, \left( \bigl| \ln|z-y| \bigr| + 1\right). 
\end{equation} 

\medskip
{\em Let us show  proposition \ref{propdecay}, postponing the proof of  lemma
  \ref{green2} to the next section}. 
 We  first need to connect the solution $v(\omega,\cdot)$ of \eqref{BL}  
to the Green function $G_\omega$. For this purpose,  we  rather consider
$$ w(\omega, y) \: = \: v(\omega,y) +  y_2 \, {\bf 1}_{\{y_2 < 0\}}(y). $$ 
which satisfies \eqref{BL2}. Note that $v$ and $w$  coincide at
$y=0$. Formally, $w$ should be equal to  
\begin{equation} \label{tildew}
 \tilde{w}(\omega,z) \: = \: \int_{\{y_2 =0\}} G_\omega(z,y) \, e_1 \, dy. 
\end{equation}
Using estimates \eqref{greenestim1}, \eqref{greenestim3}, it is standard to
show that $\tilde{w}$ is a solution of \eqref{BL2} in 
$H^1_{loc}(\overline{\Omega_{bl}})$.
Using  bound \eqref{greenestim2}, one has even
$$   \int_{\Omega^{bl} \cap \{|z_1| < 1\}} | \na \tilde{w} |^2 \: \le \: C 
\: < \:  +\infty, $$
for all $\omega \in P_{\alpha}$.

\medskip
Extending $\tilde{w}(\omega, \cdot)$ by $0$  outside
$\Omega^{bl}(\omega)$, one can show that $\omega \mapsto
\tilde{w}(\omega,\cdot)$ is measurable from $P_{\alpha}$ to
$H^1_{loc}(\R^2)$ (see the appendix  for details). 
 Moreover, thanks to \eqref{stationG}, $\tilde{w}$ satisfies the stationarity
  relation $\tilde{w}(\tau_h(\omega),y) = \tilde{w}(\omega,
  (y_1+h,y_2))$. Finally, using that $w$ and $\tilde{w}$ both satisfy
 \eqref{BL2}, a simple energy estimate on the difference leads to  
$$ \E \int_{\Omega^{bl} \cap \{|z_1| < 1\}} \left| \na \left( \tilde{w} - w
  \right) \right|^2 \: = \: 0 $$
which shows that $w = \tilde{w}$ almost surely.

\medskip
It remains to estimate the difference
$$ v(\omega_1,0,0) - v(\omega_2,0,0) \: = \:   \int_{\{y_2 =0\}}
\left( G_{\omega_1}((0,0),y) -  G_{\omega_2}((0,0),y) \right) \, e_1,   $$
for  every $\omega_1$, $\omega_2$ in $P_\alpha$
  which coincide over $[-n,n]$. This integral is bounded  by
\begin{align*}  
I_1 + I_2 \: & :=  \: \int_{y_2=0, |y_1| > n} \left| G_{\omega_1} -
G_{\omega_2}\right|((0,0),y) \, dy \\
&  + \: \int_{y_2=0, |y_1| \le  n} \left| G_{\omega_1} -
G_{\omega_2}\right|((0,0),y) \, dy
\end{align*}
The use of \eqref{greenestim1} gives
$$ I_1 \: \le \: C \, \int_{y_2=0, |y_1| > n} \frac{1}{|y_1|^{2\tau}} \,
 dy_1 \: \le \: \frac{C}{n^{2\tau-1}}. $$
where $C$, which  depends {\it a priori}  on  $\| \omega_i
\|_{C^{2,\alpha}}$, can be chosen uniformly over  $P_\alpha$, as all
$C^{2,\alpha}$ norms are bounded by $K_{\alpha}$. 
 To bound the second term, we first assume that $\omega_2 > 
 \omega_1$ for $|y_1| >  n$, which is always possible up to introduce an
 intermediate third
 boundary. Hence, $\Omega^{bl}(\omega_2) \subset \Omega^{bl}(\omega_1)$. To
 lighten notations, we introduce 
$$ \Omega^{bl}_{1,2} \: := \:
 \Omega^{bl}(\omega_1)\setminus\Omega^{bl}(\omega_2), \quad \Gamma_{1,2} \:
 := \:  \pa  \Omega^{bl}(\omega_2)\setminus \pa \Omega^{bl}(\omega_1), $$
as well as  
$$  \tilde{P}(y) \: := \: 
P_{\omega_2}\bigl((0,0),(y_1,\omega_2(y_1))\bigr), \quad y \in \Omega^{bl}_{1,2}$$
which defines a continuous extension of  $P_{\omega_2}((0,0), \cdot)$
 outside $\Omega^{bl}(\omega_2)$. Finally, we define the vector fields
\begin{align*}
& U(y) \:  := \: \left(G_{\omega_1} - G_{\omega_2}\right)((0,0),y),  \quad
   Q(y) \: := \:  \left(P_{\omega_1}-P_{\omega_2}\right)((0,0),y), \quad y \in
\Omega^{bl}(\omega_2), \\
&  U(y)  \:  := \:  G_{\omega_1}((0,0),y), \hspace{1.9cm}   Q(y) \: := \:
  P_{\omega_1}((0,0),y) - \tilde{P}(y),  \quad y \in 
\Omega^{bl}_{1,2}.
\end{align*} 
They satisfy
\begin{equation}
\left\{
\begin{aligned}
-\Delta U + \na Q & = 0,  \quad  y \in \Omega^{bl}(\omega_2), \\
-\Delta U + \na Q & = -\na \tilde{P},    \quad y \in
\Omega^{bl}_{1,2}, \\
\div U & = 0, \quad y \in \Omega^{bl}(\omega_1), \\
U & = 0,  \quad y \in \pa \Omega^{bl}(\omega_1), \\
\left[ U \right]\vert_{\Gamma_{1,2}}&  = 0,  \quad \left[ \pa_n U -
  Q \otimes n \right]\vert_{\Gamma_{1,2}} =
-\pa_n G_{\omega_2}((0,0),y)\vert_{\Gamma_{1,2} }.
\end{aligned}
\right.
\end{equation}
  A direct energy estimate yields 
\begin{align*}
& \int_{\Omega^{bl}(\omega_1)} \left| \na U \right|^2 \: \le \:
 \int_{\Gamma_{1,2}} |\pa_n G_{\omega_2}((0,0),y)| \, |U| \: + \: 
\int_{\Omega^{bl}_{1,2}} |\na \tilde{P}| \,  |U|  \\
&  \: \le \:  \Bigl( \int_{\Gamma_{1,2}} |\pa_n
 G_{\omega_2}((0,0),y)|^2 \Bigr)^{1/2}  \Bigl(  \int_{\Gamma_{1,2}}  |U|^2
 \Bigr)^{1/2} 
 \: + \:  \Bigl( \int_{\Omega^{bl}_{1,2}} |\na \tilde{P}|^2
 \Bigr)^{1/2} \Bigl(  \int_{\Omega^{bl}_{1,2}}  |U|^2
 \Bigr)^{1/2} \\
&  \: \le \: C \, \biggl(  \Bigl( \int_{\Gamma_{1,2}} |\pa_n
 G_{\omega_2}((0,0),y)|^2 \Bigr)^{1/2} + \Bigl(\int_{\Omega^{bl}_{1,2}} | \na \tilde{P}
 |^2 \Bigr)^{1/2} \biggr) \Bigl( \int_{\Omega^{bl}_{1,2}}
 | \pa_{y_2} U |^2 \Bigr)^{1/2}. 
\end{align*}
Note that all $y$ in  both $\Gamma_{1,2}$ and $\Omega^{bl}_{1,2}$
 satisfy $|y_1| > n$. Using \eqref{greenestim1}, we end up with 
$$\int_{\Omega^{bl}(\omega_1)} \left| \na U \right|^2 \: \le \: C \, n^{1-4\tau}$$
Back to $I_2$, we obtain
\begin{equation*}
 | I_2 | \:  \le \:  \sqrt{2 n} \, \left( \int_{|y_1| \le n, y_2 = 0} \,
   |U|^2  \right)^{1/2} \:
  \le \: C \, \sqrt{n} \,  \left( \int_{\Omega^{bl}(\omega_1)} \, 
   |\pa_{y_2} U|^2  \right)^{1/2} \: \le \: \frac{C}{n^{2\tau-1}}.  
\end{equation*}
This ends the proof of proposition \eqref{propdecay}.

\section{Green function estimates}
{\em This section is devoted to the proof of lemma
  \ref{green2}}, that is sharp estimates on the Green function
  $(G_\omega,P_\omega)$ for the Stokes operator above the humped boundary
  $y_2 = \omega(y_1)$, where $\omega$ belongs to $C^{2,\alpha}$. A
  fundamental remark is that the Green function satisfies the scaling 
\begin{equation} \label{scaling}
\forall \eps > 0, \quad G_{\omega^\eps}(\eps z, \eps y) \: = \:
G_\omega(z,y), \quad \omega^\eps(x_1) = \eps \omega(x_1/\eps). 
\end{equation}
We want estimates \eqref{greenestim1} to hold for $|z-y|$ large,
that is for $\eps \: :=  |z-y|^{-1}$ small. By relation \eqref{scaling},
to establish such estimates amounts to get local estimates for the
Green function $G_{\omega^\eps}$. Thus, this is again a homogenization
problem: more precisely, we must show that the oscillations of the boundary
at frequency $\eps^{-1}$ do not affect too much the estimates on
$G_{\omega^\eps}$,  so that it
behaves as the Green function for a half-plane. A very close problem has
been considered in the papers \cite{Avellaneda:1987, Avellaneda:1991}
 by Avellaneda and Lin, namely the
derivation of local estimates for elliptic systems $\div\left(A(x/\eps) \na
\, \cdot  \, \right)$, in which $A$ is a positive definite matrix with
periodic coefficients. Our reasoning follows  these papers.

 For all $x \in \R^2$, $r > 0$, we will denote $D(x,r)$
the disk of center $x$ and radius $r$, and 
\begin{equation*}
 D^\eps(x,r) \: := \:
D(x,r) \cap \{x_2 > \omega^\eps(x_1)\}, \quad 
\Gamma^\eps(x,r)\: := \: D(x,r) \cap \{x_2 = \omega^\eps(x_1)\}.  
\end{equation*}
An important property is that for all $0 < r < 1$, 
\begin{equation} \label{airedeps}
| D^\eps(x,r) | \: \ge \: \eta \, r^2, 
\end{equation} 
for some $\eta > 0$ independent of $\eps$. More precisely, $\eta$ only
involves the Lipschitz norm of $\omega^\eps$, wich is bounded uniformly in
$\eps$. This will be used implicitly throughout the sequel. 

The core of the proof is to derive elliptic estimates uniform with
respect to $\eps$ on the following Stokes problem:
\begin{equation} \label{stokes}
\left\{
\begin{aligned}
-\Delta u  + \na p \: = \: \div  f, \quad & x \in D^\eps(x_0,1)\\
\div u = 0, \quad & x \in D^\eps(x_0,1),\\
u = 0, \quad  & x \in \Gamma^\eps(x_0,1), 
\end{aligned}
\right.
\end{equation}
where $x_0 \in \R^2$. More precisely, there are  two steps in the proof:
\begin{enumerate}
\item We show a $\eps$-uniform H\"older estimate on $u$: for all $f \in
  L^q$, $q>2$  and for $\mu = 1-2/q$,  
\begin{equation} \label{holder}
 \| u \|_{C^{0,\mu}(D^\eps(x_0,1/2))} \: \le \: C \left( \| f
  \|_{L^q(D^\eps(x_0,1))} + \| u \|_{L^2(D^\eps(x_0,1))} \right). 
\end{equation} 
where $C$ depends only on $\| \omega \|_{C^{1,\alpha}}$.
\item  Thanks to this H\"older estimate, we  prove \eqref{greenestim1}. 
\end{enumerate}
The two next paragraphs correspond to these steps.
\subsection{H\"older estimate}
To obtain a H\"older regularity result, a classical approach
is to use  a characterization of H\"older spaces due to Campanato (see
\cite{Giaquinta:1983}): for $\Omega$ an open connected bounded set,  
$\displaystyle v \in C^{0,\mu}(\Omega)$ iff  $v \in L^2(\Omega)$ and  
$$  \sup_{x \in \Omega, r > 0} \frac{1}{r^{2+2\mu}}\int_{\Omega(x,r)} |v -
\overline{v}_{x,r}|^2 < \infty, \quad  \Omega(x,r) := \Omega \cap D(x,r), \:\:
\overline{v}_{x,r} := \frac{1}{|\Omega(x,r)|} \int_{\Omega(x,r)} v. $$
One then tries to control such local integrals through energy
estimates. This approach 
has been successful in the study of elliptic systems, see the work of
Giaquinta and coauthors \cite{Giaquinta:1983}. It extends to the Stokes type
equations, {\it cf} article \cite{Giaquinta:1982}.
 For us, it amounts to controlling 
$$   I^\eps_{x,r} \: := \: \frac{1}{r^{2+2\mu}}\int_{D^\eps(x,r)} |u - 
\overline{u}_{x,r}|^2 <
\infty, \quad  \overline{u}_{x,r} := \frac{1}{|D^\eps(x,r)|}
\int_{D^\eps(x,r)} u $$
where $u$ is solution of \eqref{stokes}. Note that, thanks to
\eqref{airedeps} (see \cite{Giaquinta:1983}),   
$$ \| u \|_{C^{0,\mu}(D^\eps(x_0,1/2))} \: \le \: C_{x_0} 
\left( \|  u \|_{L^2(D^\eps(x_0,1/2))}  \: +   
\sup_{x \in D^\eps(x_0,1/2), r > 0} I^\eps(x,r) \right) 
 $$
with $C_{x_0}$ independent of $\eps$. 
In our case, the main  problem is to keep track of the dependence of
 $I^\eps_{x,r}$ on 
 $\eps$. It involves a  discussion of the way $\eps$ relates to $r$. 
Broadly speaking, the idea is the following: if $r$ is large compared to 
$\eps$, then the oscillations have small enough amplitude to apply the
 regularity results of the flat case. On the contrary,  when $r$ gets as
 small  or even smaller than $\eps$,
one can rescale everything by a factor $\eps$, so that oscillations of the
boundary have frequency $O(1)$, and are no longer annoying. 
Implementation of this  idea is a bit technical, and follows closely 
the work of Avellaneda and Lin. 

\medskip
We first remind a few elements of regularity theory for Stokes type
systems. Let $\Omega$ an open connected bounded set, with {\em 
Lipschitz boundary}. Then, for any $\varphi \in L^2(\Omega)$ satisfying
$\int_\Omega \varphi = 0$, the problem 
$$ \div w = \varphi, \quad w\vert_{\pa \Omega} = 0 $$
has  one solution $w$ satisfying $\| w \|_{H^1_0} \: \le \:  C \, \| \varphi
\|_{L^2(\Omega)}$, where {\em  $C$ can be taken as an increasing function of
$|\Omega|$ and  of the Lipschitz constant $K$ of the boundary},
 see \cite{Galdi:1994}. 
Thanks to this result, one has quite easily, 
\begin{equation} \label{pressure}
\| p - \int_\Omega p \|_{L^2(\Omega)} \: \le \: C \, \| \Delta u +f \|
\|_{H^{-1}(\Omega)} 
\end{equation}
where $(u,p) \in H^1(\Omega) \times L^2(\Omega)$ satisfies (in the
distributional sense) 
\begin{equation} \label{stokesbis}
 -\Delta u + \na p = f, \quad \div u = 0, \quad x \in \Omega.
\end{equation}
Again, the constant $C$ in \eqref{pressure} depends only on $|\Omega|$ and
the Lipschitz constant of the boundary. 

\medskip
We now state  the famous Cacciopoli inequality: 
\begin{lemma}
For all $0 < r < 1$,
 any solution $u \in H^1(\Omega)$ of \eqref{stokes} satisfies 
\begin{equation} \label{cacciopoli}
\| \na  u \|_{L^2(D^\eps(x,r))} \: \le \:  C \left( r^{-1} \, \| 
u \|_{L^2(D^\eps(x,2r))} \: + \: r^\mu \, \| f \|_{L^q(D^\eps(x,2r))} \right).
\end{equation}
\end{lemma}
{\em Sketch of proof:} We remind the main elements of proof.
 Let $\eta$ a smooth function with compact support in
$D(x,2r)$, with $\eta=1$ on $D(x,r)$. Note that $|\na \eta| \le C r^{-1}$. 
Multiplying \eqref{stokes} by the test function $\eta^2 u$, and integrating
by parts, one has easily
\begin{align*}
 \int_{D^\eps(x,r)}|\na u|^2 \: & \le \: \int_{D^\eps(x,2r)}\eta^2 |\na
u|^2 \:  \le \: C \, r^{-2} \, \int_{D^\eps(x,2r)} |u|^2 \: + \:
C \, \int_{D^\eps(x,2r)} |f|^2 \\
& \: + \: 
 \| p - \overline{p}_{x,2r} \|_{L^2(D^\eps(x,2r))}
\, \| \div (\eta u) \|_{L^2(D^\eps(x,2r))}.
\end{align*}
Using \eqref{pressure}, we get 
$$  \| p - \overline{p}_{x,2r} \|_{L^2(D^\eps(x,2r))} \le 
C \, \| \Delta u + \div f \|_{H^{-1}(D^\eps(x,2r))} = C \| v
\|_{H^1(D^\eps(x,2r))}  $$
 where $v \in H^1_0(D^\eps(x,2r))$ is the solution of 
$$-\Delta v + \na p =  \Delta u + \div f, \quad \div v = 0, 
\quad v\vert_{\pa D^\eps(x,2r)} = 0. $$
Note that the previous bound is uniform in $\eps$, as it only involves the
Lipschitz constant of $\omega^\eps$ which is uniformly bounded.
A simple energy estimate on $v$ gives
$$ \| \na v  \|_{L^2(D^\eps(x,2r))} \: \le \: C \left(  \| \na u
\|_{L^2(D^\eps(x,2r))} \: + \: \| f \|_{L^2(D^\eps(x,2r))} \right) $$
As $\div(\eta u) = \na \eta \cdot u$, and using H\"older inequality on $f$,
we end up with 
\begin{align*}
 \int_{D^\eps(x,r)}|\na u|^2 \:  \le \: \int_{D^\eps(x,2r)}\eta^2 |\na
u|^2 \:  & \le \: C \, r^{-2} \, \int_{D^\eps(x,2r)} |u|^2 \: + \:
C_\delta \, r^{2\mu} \,  \|f\|^2_{L^q(D^\eps(x,2r))} \\
& \: + \: 
 \delta \| \na u \|^2_{L^2(D^\eps(x,2r))},
\end{align*}
where $\delta > 0$ is arbitrary small. We conclude as in \cite[Theorem 1.1,
page 180]{Giaquinta:1982}. 

\medskip
Inequality of type \eqref{cacciopoli} has been used by Giaquinta and Modica
in the study of elliptic regularity. In the context of Stokes type system,
they obtain a local estimate, see \cite{Giaquinta:1982}:
\begin{theorem} \label{modica}
Let $\Omega$ of class $C^1$, and 
 $(u,p,f) \in H^1(\Omega) \times L^2(\Omega) \times L^q(\Omega)$, $q > 2$, 
  satisfying
$$  -\Delta u + \na p = \div f, \: \div u = 0, \quad x \in \Omega(x_0,1),
 \quad u\vert_{\pa \Omega \cap D(x_0,1)} = 0. $$
 Then, $u \in  C^{0,\mu}(\Omega(x_0,1/2))$ for $\mu=1-2/q$, and  
\begin{equation} \label{modicaeq}
\| u \|_{C^{0,\mu}(\Omega(x_0,1/2))} \: \le \: C 
\left( \| u \|_{L^2(\Omega(x_0,1))} \: + \: \| f \|_{L^q(\Omega(x_0,1))}
 \right). 
 \end{equation}
\end{theorem}
Unfortunately, we cannot use this theorem assuch. Indeed, the constant $C$
in the last regularity estimate involves the modulus of continuity of $\na
\gamma$, where $x_2 = \gamma(x_1)$ describes the boundary. In our case
$\gamma = \omega^\eps$,  such modulus of continuity is not uniformly
bounded in $\eps$. We must proceed in several steps to control
the local integrals $I^\eps_{x,r}$. Note that theorem \ref{modica} implies 
estimate  \eqref{holder}  when $D^\eps(x_0,1)$
is far from the boundary. Thus, we
can restrict ourselves to a case in which  $x_0$ is close to the
oscillating boundary, for instance  belongs to the axis $x_2 = 0$. 

\begin{lemma} \label{lemme1}
For all  $\theta$ small enough, there exists $\eps_0 > 0$ 
such that for all $\eps < \eps_0$,  and 
for all solutions of \eqref{stokes} satisfying 
 $\displaystyle \quad    \| u^\eps \|_{L^2(D^\eps(x_0,1/4))} \le 1$,  
$ \displaystyle \:  \| f \|_{L^q(D^\eps(x_0,1/4))} \le \eps_0$, one has  
$$\| u^\eps \|_{L^2(D^\eps(x_0,\theta))}
 \, \le \,  \theta^{\mu+1}.  $$
\end{lemma}
{\em Proof of the lemma:} Suppose that the result does not hold. Then one
can find $\theta$ arbitrary small, and 
sequences  $u^{\eps_j}$, $f^j$ satisfying  
$$ \| u^{\eps_j} \|_{L^2(D^{\eps_j}(x_0,1/4))} \le 1, \quad 
\| f^j \|_{L^2(D^{\eps_j}(x_0,1/4))} \xrightarrow[j \rightarrow +\infty]{}
0, \quad \| u^{\eps_j} 
 \|_{L^2(D^{\eps_j}(x_0,\theta))} > \theta^{\mu+1}.$$
 One can extend all the 
$u^\eps_j$, $f^j$ by $0$ outside $D^{\eps_j}(x_0,1/4)$
so that all these functions are defined on the fixed domain
$D(x_0,1/4)$. From the $L^2$ bound on $u^{\eps_j}$, up to extract a 
subsequence, we get 
$$   u^{\eps_j} \mbox{ converges weakly to some } u \mbox{ in } 
 L^2(D(x_0,1/4)) $$  
and by Cacciopoli inequality \eqref{cacciopoli}, 
$$  u^{\eps_j} \mbox{ converges weakly  to  } u  \mbox{ in }
H^1(D(x_0,1/8)),  \mbox{ and strongly  in }
L^2(D(x_0,1/8)). $$
One can then take the limit in \eqref{stokes}, which yields 
$$ -\Delta u + \na p = 0, \quad \div u = 0, \quad \mbox{in } 
 D(x_0,1/8) \cap \{x_2 >
0\}, \quad u\vert_{D(x_0,1/8) \cap \{ x_2 = 0 \}} =0. $$ 
As the upper half plane is a regular domain, we can apply theorem
\ref{modica}, so that for all $\tilde{\mu} > \mu$, for all $\theta$,  
\begin{align*}
 \| u \|_{L^2(D(x_0,\theta) \cap  \{x_2 > 0\})} \: & \le  \: 2 \pi \, \| u \|_{C^{0,\tilde{\mu}}(D(x_0,\theta) \cap  \{x_2 > 0\})} 
\, \theta^{\tilde{\mu}+1} \\
&  \le \: C \,\| u \|_{L^2(D(x_0,1/8)\cap  \{x_2 > 0\})}  \,
\theta^{\tilde{\mu}+1} \: \le \: C \, \theta^{\tilde{\mu}+1}.  
\end{align*}
For $\theta$ small enough, it  contradicts the lower bound on
 $\| u^{\eps_j} \|_{L^2(D^{\eps_j}(x_0,\theta))}$.   

\medskip
We fix $\theta$, $\eps_0$ as in 
lemma \ref{lemme1}. We then state 
\begin{lemma} \label{lemme2}
For all $\eps, k$ satisfying $\eps/\theta^k \le \eps_0$ ($k \ge 1$), 
\begin{equation} \label{localk}
\int_{D^\eps(0,\theta^k)} |u^\eps |^2 \:
 \le \: \theta^{2k\mu+2}   \, \left( \| u^\eps
 \|_{L^2(D^\eps(0,1/4))} \: + \: 
\frac{1}{\eps_0} \, \| f \|_{L^2(D^\eps(0,1/4))} \right)^2
\end{equation} 
\end{lemma}  
{\em Proof of the lemma:} 
The lemma is deduced from an induction argument on $k$. For $k=1$,
the bound \eqref{localk} is given by lemma \ref{lemme1}.
 Assume now that this bound holds for $k \ge 1$. Up to
 a horizontal translation, we can assume that $x_0 =0$. 
Then, we set 
$$M  := \| u^\eps \|_{L^2(D^\eps(0,1/4))} \: + \: 
\frac{1}{\eps_0} \, \| f \|_{L^2(D^\eps(0,1/4))},$$ 
and introduce the rescaled functions 
$$ v \: := \: \theta^{-k\mu} M^{-1}  u(\theta^k x), \quad g \: = \:
 \theta^{k-k\mu} M^{-1} f(\theta^k x).   $$
They satisfy 
$$ -\Delta v + \na q = \div g, \quad \div v = 0, \: x \in
D^\eps(0,\theta^{-k}/4), \quad v\vert_{\Gamma^\eps(0,\theta^{-k}/4)}.$$
Moreover, one has $\displaystyle \| f \|_{L^q(D^\eps(0,1/4))} \, \le \,
\eps_0$, and by the induction assumption 
$$\displaystyle \| v \|_{L^2(D^\eps(0,1/4))} \, \le \,  1.$$ 
Applying  lemma \ref{lemme1} to $v$
and $g$ yields the result. 

\medskip
We can now finish the proof of estimate \eqref{holder}. 
Let $x \in D^\eps(x_0,1/2)$. We need to bound $I^\eps(x,r)$, for $r > 0$.
 There are two cases to handle:
\begin{itemize}
\item {\em The distance between  $x$ and  the boundary $\{ x_2 =
  \omega^\eps(x_1)\}$  satisfies $\delta^\eps(x) \ge \frac{\eps}{\eps_0}$.} 
\end{itemize}
Up to take a smaller $\eps_0$, we can suppose that $\frac{\eps}{\eps_0} >
\eps$, which implies that there exists $x'_0$ on the axis $\{x_2=0\}$ with
$|x - x'_0| \le \delta^\eps(x)$. 
By lemma \ref{lemme2},  for all $\eps/\eps_0 \le r \le 1/12$, 
\begin{equation} \label{deps3r}
 \int_{D^\eps(x'_0, 3r)} | u^\eps |^2 \: \le \: C \, r^{2\mu+2}   \left( \|
 u^\eps \|_{L^2(D^\eps(x'_0,1/4))} \: +  \: \| f \|_{L^q(D^\eps(x'_0,1/4))}
 \right)^2.  
\end{equation}
\medskip
If $r > \delta^\eps(x) /2$, $\: D^\eps(x,r) \subset D^\eps(x'_0,3r)$, and
the previous line implies  
$$  \int_{D^\eps(x, r)} | u^\eps |^2 \: \le \: C \, r^{2\mu+2}   \left( \|
u^\eps \|_{L^2(D^\eps(x'_0,1/4))} \: +  \: \| f \|_{L^q(D^\eps(x'_0,1/4))}
\right)^2. $$  
On the contrary, if $r \le \delta^\eps(x) /2$, then $\: D^\eps(x,r) \, = \,
D(x,r)$ (it does not intersect the boundary). A simple rescaling of
\eqref{modicaeq}  yields  
$$ \| u \|_{C^{0,\mu}(D(x,\delta^\eps(x)/2))} 
\: \le \:  C \left( \delta^{\eps}(x)^{-1-\mu} \| u^\eps \|_{L^2(D(x,
  \delta^\eps(x)))} 
 \: + \:  \| f \|_{L^q(D(x, \delta^\eps(x)))} \right).   $$
Thus, 
$$ \int_{D^\eps(x,r)} | u^\eps|^2  \:  \le \:  C \, r^{2\mu+2} \, \left(
\delta^\eps(x)^{-1-\mu} \| u^\eps \|_{L^2(D(x,  \delta^\eps(x)))} \; + \:
\| f \|_{L^q(D(x, \delta^\eps(x)))} \right)^2.   $$ 
Now, by lemma \ref{lemme2}, as $\delta^\eps(x) \ge \eps/\eps_0$,  
\begin{align*}
 \| u^\eps \|_{L^2(D(x,  \delta^\eps(x)))} 
  \: & \le \: \| u^\eps \|_{L^2(D(x'_0,  2 \delta^\eps(x)))}     \\
     & \le \:
  C \, \delta^\eps(x)^{\mu+1} \, \left( \| u^\eps
  \|_{L^2(D^\eps(x'_0,1/4))} \: +  \: \| f \|_{L^q(D^\eps(x'_0,1/4))}
  \right). 
\end{align*}
Using the two last inequalities, we end up again with
$$  \int_{D^\eps(x, r)} | u^\eps |^2 \: \le \: C \, r^{2\mu+2}   \left( \|
u^\eps \|_{L^2(D^\eps(x'_0,1/4))} \: +  \: \| f \|_{L^q(D^\eps(x'_0,1/4))}
\right)^2, $$  
which in turn clearly implies 
$$  \int_{D^\eps(x, r)} | u^\eps - \overline{u^\eps}_{x,r} |^2 \: \le \: C
  \, r^{2\mu+2}  
  \left( \| u^\eps \|_{L^2(D^\eps(x'_0,1/4))} \: +  \: \| f
  \|_{L^q(D^\eps(x'_0,1/4))} \right)^2, $$  
As $D^\eps(x'_0,1/4)) \: \subset \: D^\eps(x_0,1)$, this gives the required
  estimate.  
\begin{itemize}
\item {\em The distance between  $x$ and  the boundary $\{ x_2 =
  \omega^\eps(x_1)\}$  satisfies $\delta^\eps(x) < \frac{\eps}{\eps_0}$.} 
\end{itemize}
This time, there exists $x'_0$ on the axis $\{ x_2 = 0 \}$ such that 
$|x - x'_0| \le \delta^\eps(x) + \eps \, \le \, 2\eps/\eps_0$. Again, for
all $\: \eps/\eps_0 \le r \le 1/12$, $\:  D(x,r)\subset D(x'_0,3r) \:$ and  
 \eqref{deps3r} implies
$$ \int_{D^\eps(x,r)} | u^\eps |^ 2 \: \le \:  C \, r^{2\mu+2}   \left( \|
u^\eps \|_{L^2(D^\eps(x'_0,1/4))} \: +  \: \| f \|_{L^q(D^\eps(x'_0,1/4))}
\right)^2. $$  
It remains to handle the case $r < \eps/\eps_0$. Up to a horizontal
translation, we can assume that $x'_0=0$. We  introduce the rescaled
functions  
 $$ v \: =  \: \left(\frac{\eps}{\eps_0}\right)^{-\mu}
u^\eps\left(\frac{\eps}{\eps_0} x\right), \quad g \: =  \:
\left(\frac{\eps}{\eps_0}\right)^{1-\mu} f\left(\frac{\eps}{\eps_0}
x\right).$$  
They satisfy  in particular
\begin{equation}
-\Delta v + \na q = \div g, \quad \div v = 0, \: x \in D^{\eps_0}(0,1),
 \quad v\vert_{\Gamma^{\eps_0}(0,1)} = 0. 
\end{equation}

It is a Stokes type system set in a domain independent of the small
parameter $\eps$.
Hence, we can apply theorem \ref{modica}, which yields: for all $r < 1$,
$$ \int_{D^{\eps_0}(x,r)} | v - \overline{v}_{x,r} |^2 \: \le \: C \, \| v
 \|_{C^{0,\mu}(D^{\eps_0}(x,r))} \, r^{2\mu+2} 
 \: \le \: C \, r^{2\mu+2} \, \left( \| v \|_{L^2(D^{\eps_0}(0,2))}  \: + \: 
\| g \|_{L^2(D^{\eps_0}(0,2))} \right)^2. $$
Back to the original unknowns $u^\eps$, $f$, we obtain the control of
 $I^\eps_{x,r}$ for $r \le \eps/\eps_0$. This ends the proof.  

\subsection{Bounds on the Green function}
From the above H\"older estimate, we can deduce the estimate
\eqref{greenestim1}. The arguments are again adapted from  article
\cite[pages 819,829-831]{Avellaneda:1987}. 
 For the sake of completeness, we describe the ideas at play.  
The first step is to establish the following bound: for all $x,x'$ in $\{
x_2 > \omega^\eps(x_1)\}$,   
\begin{equation} \label{logG}
 | G_{\omega^\eps}(x,x')  | \: \le \: C \, \left( \bigl| \ln |x-x'| \bigr|
 \: + \: 1
   \right).  
 \end{equation}
where $C$ only involves $\| \omega \|_{C^{1,\alpha}}$. Note that it implies
\eqref{greenestim3}. 
To lighten the notations, we drop the suffix $\omega$, denoting $G^\eps$,
$G$ instead of $G_{\omega^\eps}, G_\omega$. Let us introduce
the  
Green function $\tilde{G}^\eps(x,t,x',t')$ for the Stokes operator over $\{
x_2 > \omega^\eps(x_1)\} \times \T$. Namely, it satisfies for all $(x,t)
\in \{ x_2 > \omega^\eps(x_1)\} \times \T$ 
\begin{equation}
\left\{
\begin{aligned}
-\Delta  \tilde{G}^\eps(x,t,\cdot) +  \na \tilde{P}^\eps(x,t,\cdot) \: = \:
 \delta_{x,t} I_3 & 
 \quad \mbox{ in }  \{ x_2 > \omega^\eps(x_1)\} \times \T, \\
  \div \tilde{G}^\eps(x,t,\cdot) = 0 & \quad \mbox{ in } \{ x_2 >
 \omega^\eps(x_1)\}  \times \T, \\  
 \tilde{G}^\eps(x,t, \cdot) = 0 & \quad  \mbox{ on }  \{ x_2 =
 \omega^\eps(x_1)\}   \times \T.  
\end{aligned}
\right.
\end{equation}
One has easily that 
$$ G^\eps(x,x') \: = \:  \int_{\T} \left( \tilde{G}^\eps_1(x,0,x',t'),
\tilde{G}^\eps_2(x,0,x',t') \right) \, dt'. $$  
The point is to show that  
$$ |\tilde{G}^\eps(x,t,x',t')| \: \le  \: C \, \frac{1}{ |x-x'| + |t-t'|}. $$
The estimate \eqref{logG} is then obtained by integration with respect to
$t'$,  at $t=0$. 
Such bound on $\tilde{G}^\eps$ will be deduced from a repeated use of  the
H\"older estimate \eqref{holder}. Note that  this estimate extends without
difficulty to  similar Stokes problems in  dimension $n \ge 2$, with  $q >
n$  and $\mu = 1 - n/q$. In particular, it holds  when the domain is $\{
x_2 > \omega^\eps(x_1)\}  \times \T$. 

\medskip
Namely, let $\tilde{x} := (x,t)$ and $\tilde{x}' := (x',t')$. Let $r :=
 |\tilde{x}  - \tilde{x}'|$, and
 $f \in C^\infty_c(D^\eps(\tilde{x}',r/3))$. We consider the quantity 
$$u^\eps(\tilde{x}) \: =  \: \int_{D^\eps(\tilde{x}',r/3)}
\tilde{G}^\eps(\tilde{x}, \tilde{z}) \, f(\tilde{z}) \, d\tilde{z}$$   
The field $u^\eps$ satisfies a Stokes equation with source term $f$  over
$\{ x_2 > \omega^\eps(x_1)\}  \times \T$, with a Dirichlet boundary
condition. We  therefore apply the estimate \eqref{holder} to
$u^\eps$. Properly rescaled, it yields 
$$  | u^\eps(\tilde{x}) | \: \le \: C \, \left( r^{-3}
\int_{D^\eps(\tilde{x},r/3)} |u^\eps|^2 \right)^{1/2}, $$ 
where we used the fact that $f$ vanishes over
$D^\eps(\tilde{x},r/3)$. Thus, we get   
\begin{align*}
& \left| \int_{D^\eps(\tilde{x}',r/3)} \tilde{G}^\eps(\tilde{x}, \tilde{z})
  \, f(\tilde{z}) \, d\tilde{z} \right| \:   \le \:  C \,  \left( r^{-3}
  \int_{D^\eps(\tilde{x},r/3)} |u^\eps|^2 \right)^{1/2}  \\ 
&  \le \: C' \,  \left( r^{-3} \int_{D^\eps(\tilde{x},r/3)} |u^\eps|^6
  \right)^{1/6} \: \le \:  
C' \, r^{-1/2} \,  \left( \int_{D^\eps(\tilde{x},r/3)} |u^\eps|^6
\right)^{1/6} \\ 
&  \le \: C'' \, r^{-1/2} \left(  \int_{\{ x_2 > \omega^\eps(x_1)\}  \times \T}
 | \na u^\eps |^2 \right)^{1/2} \: \le \:  C'' \, r^{1/2}  \| f \|_{L^2}^{1/2} 
\end{align*}
where the last two inequalities come respectively from the Sobolev
imbedding theorem (note that $u^\eps$ is zero at the boundary so that the
imbedding does not involve lower order terms), and from the standard energy
estimate on the Stokes system.  
By duality, we infer that 
$$ \left( r^{-3} \int_{D^\eps(\tilde{x}', r/3)} |\tilde{G}^\eps(\tilde{x},
\tilde{z})|^2  \, d\tilde{z} \right)^{1/2} \: \le \: C \, r^{-1}. $$ 
Using that $\tilde{G}^\eps(\tilde{x}, \cdot)$ satisfies a homogenenous
Stokes system over $D(\tilde{x}',r/3)$, one more application of 
\eqref{holder} leads to 
$$ \tilde{G}^\eps(\tilde{x}, \tilde{y}) \: \le \: C \, \left( r^{-3} 
\int_{D^\eps(\tilde{x}', r/3)} |\tilde{G}^\eps(\tilde{x}, \tilde{z})|^2  \,
 d\tilde{z} \right)^{1/2} \: \le \: C \, r^{-1}. $$

\medskip
Inequality \eqref{logG} at hand, we can derive the final estimate on
$G$. Let  $x, x' \in \{ x_2 > \omega^\eps(x_1)  \}$. Set this time 
$r := |x - x'|$. For all $\bar{x}$   such that 
$|\bar{x}-x| < 2r$,  \eqref{logG} implies 
$$ \left| G^\eps(\bar{x},x') \right| \: \le \: C \left( \bigl| 
\ln |\bar{x}-x'| \bigr| +
1 \right)\: \le \: C' \left( \bigl| \ln |x-x'| \bigr|  +1 \right). $$
Applying \eqref{holder} to the function $G^\eps(\cdot,x')$, we get for any
$\tau \in (0,1)$  
\begin{align*}
 \left| G^\eps(x,x') \right| \: & \le \:   \, \delta^\eps(x)^\tau 
\, \| G^\eps(\cdot ,x') \|_{C^{0,\tau}(D(x,r/3))} \\
&  \le \: C_\tau \,
\delta^\eps(x)^\tau \, r^{-1-\tau} \, \|  
 G^\eps(\cdot ,x') \|_{L^2(D(x,2r/3))} 
\end{align*} 
that leads to 
$$  \left| G^\eps(x,x') \right| \:  \le \: C_\tau \, \left( \bigl| 
\ln |x-x'| \bigr| +1 \right) \,  
\frac{\delta^\eps(x)^\tau}{|x-x'|^\tau}. $$ 
Now reversing the roles of  $x$ and $x'$, we obtain  
$$ |G^\eps(x,x')| \: \le \: C_\tau \, \left( \bigl| \ln |x-x'| \bigr|
 +1 \right)
\frac{\delta^\eps(x)^\tau \, \delta^\eps(x')^\tau}{|x-x'|^{2\tau}}. $$ 
Using the scaling relation  \eqref{scaling}, we get for all $\: y,z \in
\{ y_2 > \omega(y_1) \}, \:$ for $\: \eps := |y-z|^{-1}$,
$$  G(z,y) = G^\eps\left(\eps z, \eps y \right) \: \le \: 
C_\tau  \left( \bigl| \ln \left(\eps |z-y| \right) \bigr| + 1 \right)
\frac{\delta^\eps(\eps z)^\tau \, \delta^\eps(\eps y)^\tau}{|\eps
  (y-z)|^{2\tau}}  = C_\tau 
\frac{\delta(z)^\tau \, \delta(y)^\tau}{|y-z|^{2\tau}}.  $$
Using  classical local regularity results  for the Stokes equation in a
$C^{2,\alpha}$ domain (see \cite[Theorem 1.3, page 198]{Giaquinta:1982},
 which extends theorem \ref{modica}):
for $|z-y| \ge 1$, 
\begin{align*} 
\sum_{|\beta| \le 2} |\pa^\beta_y G(z,y)| \: + \: |\na_y
 P(z,y)| \:  & \le  \: C \, \| G(z,\cdot) \|_{L^2(D(y,1/2))} \\
& \le  \: C \,  \frac{ \delta(z)^\tau \, (1+ \delta(y))^\tau}{|z-y|^{2\tau}}
\end{align*}
that is exactly estimate \eqref{greenestim1}. 

\section{Comments}
\subsection{Well-posedness of the boundary layer system}
As mentioned several times in this paper, 
the well-posedness of system \eqref{BL} is not known 
without a structural assumption on  $\omega$, like
periodicity or stationarity. We stress however that thanks to our estimates
on the Green function $G_\omega$, the representation formula \eqref{tildew}
defines a solution of \eqref{BL} for any $C^{2,\alpha}$ boundary, {\em cf}
the fourth section. Hence, the open  issue is rather to find the appropriate
functional space for uniqueness. 

\medskip
Such difficulty does not arise when the Stokes operator is replaced by  
the Laplacian, or more generally by a scalar elliptic operator. Hence,
 one can show well-posedness in $L^\infty$ of 
\begin{equation} \label{BLbis}
 -\Delta v = 0  \: \mbox{ in } \Omega^{bl},  
 \quad v(y) = \omega(y_1)  \: \mbox{ on  } \pa \Omega^{bl} 
\end{equation}  
if the function $\omega$ is bounded and Lipschitz.
For the existence part, one may consider, for all $n \ge 1$, the solution  $v^n$  of 
$$ -\Delta v^n = 0  \mbox{ in } \Omega^{bl} \cap D(0,n),  
 \quad v^n(y) = \omega(y_1)  \: \mbox{ on  } \pa \left( \Omega^{bl} \cap
 D(0,n)\right).  $$ 
By the maximum principle, $\| v^n \|_{L^\infty} \le \| \omega
 \|_{L^{\infty}}$, so that up to a subsequence it converges to some $v$ in 
  $L^\infty$ weak*. Straightforwardly,  $v$ satisfies \eqref{BLbis}. 

For the uniqueness part, let $v \in L^\infty$ satisfying 
\begin{equation*} 
 -\Delta v = 0  \: \mbox{ in } \Omega^{bl},  
 \quad v(y) = 0  \: \mbox{ on  } \pa \Omega^{bl} 
\end{equation*}  
Let us show that $v=0$. As 
 we do not know the behavior of $v$ at
 infinity, it does not follow directly from the maximum principle.   
  In the case of a Lipschitz boundary $\omega$, we can conclude to the
uniqueness in the following way. By the change of variables
  $y := \left(y_1, \, y_2 - \omega(y_1)\right)$, the previous
  equation becomes 
$$ \div(A(y) \na v) = 0, \: y_2 > 0, \quad v\vert_{y_2 =0} = 0,  $$
for some elliptic matrix $A = (a_{ij})$ with bounded coefficients.
 We extend $A$ and $v$ to $\{y_2 < 0\}$ by  the formulas, $\: v(y_1, y_2)  := -v(y_1, -y_2), \:$ 
 \begin{align*}
 a_{in}(y_1, y_2) &  := -a_{in}(y_1, -y_2), \quad
 a_{nj}(y_1, y_2)   := -a_{nj}(y_1, -y_2), \quad i,j  \neq n, \\
a_{ij}(y_1, y_2)    & := a_{ij}(y_1, -y_2)\quad  \mbox{ otherwise. } 
\end{align*}
In this way, we get $ \div(A(y) \na v) = 0$ on all $\R^2$.
Harnack's inequality for elliptic equations (see \cite{Gilbarg:2001})
 leads to $$ \sup_{|y| < R} (M + v) \: \le \: C \, \inf_{|y| < R} (M + v) $$
 for any $R > 0$ and $M$ such that $M+u \ge 0$. With $M = max(0,-\inf v)$
 and $R$ going to infinity, we obtain that $v = 0$. 

\subsection{Decay of correlations}
A key element in the paper is the estimate \eqref{greenestim1} on the Green
 function for the Stokes operator above an oscillating boundary. This
 estimate relies itself on the H\"older regularity result
 \eqref{holder}. In fact, with some more calculations in the same spirit,
 one could show a refined bound:  for a $C^{1,\alpha}$ boundary $\omega$
 and a source term  
 $f\in C^{0,\mu}(D^\eps(0,1)) \:$  ($\alpha, \mu >0)$, the solution
 $u^{\eps}$ of \eqref{stokes}satisfies 
 $$ \| \na u^\eps \|_{L^\infty(D^\eps(0,1/2))} \: \le \: C \left( \| u
 \|_{L^2(D^\eps(0,1))} \: + \:  
 \| f \|_{C^{0,\mu}(D^\eps(0,1))}  \right), $$
with $C$ independent of $\eps$. This yields an optimal
 bound on the Green function, that is 
\begin{equation} 
\sum_{|\beta| \le 2} |\pa^\beta_y G_\omega(z,y)| \: + \: |\na_y
 P_\omega(z,y)| \:  \le 
 \: C \,  \frac{ \delta(z) \, (1+ \delta(y))}{|z-y|^2}.
\end{equation}
Estimate \eqref{decaybound} can in turn be improved as 
$$|v(\omega_1,0,0) - v(\omega_2,0,0)| \: \le \: C \, n^{-1}.$$
Such bound is far easier to prove when \eqref{BL} is replaced by the scalar
system \eqref{BLbis}. In this case, one may use a representation  in terms
of the standard two-dimensional  Brownian motion $ B(m,t) = \left(
B_1(m,t),  B_2(m,t) \right)$. If we denote    $(M, {\cal M},\mu)$  the
probability space on which this Brownian motion is defined,    
$$ v(\omega,0,0) \: = \:   
\int_M \left(-\omega,0\right)\bigl(B_1(m, \tau(m)) \bigr) \, d\mu $$
where $\tau$ is the exit time from $\Omega^{bl}(\omega)$ (see
\cite{Varadhan:1968}). We now want to bound
$$ v(\omega_1,0,0) - v(\omega_2,0,0) = \int_M  \, \Bigl(
\omega_1\bigl(B_1(m, \tau_1(m)) \bigr) \, - \omega_2\bigl(B_1(m, \tau_2(m))
\bigr) \,  \Bigr) d\mu$$ 
where $\tau_i$ is the exit time from $\Omega^{bl}(\omega_i)$.  
We remind that  $\omega_1 = \omega_2$ over $[-n,n]$, so  the exit times are
the same for brownian particles leaving in the region $y_1 \in
[-n,n]$. Hence, 
\begin{equation*}
\begin{aligned}
  \left| v(\omega_1,0,0) -  v(\omega_2,0,0) \right| \: & \le \: 2 \, 
  \max_{i=1,2} \: \| \omega_i\|_{L^\infty}   \: \mu\left( |B_1(\cdot,
  \tau_i)| \ge n \right)  \\ 
& \le \:  2 \, \Bigl( 
 \mu\left(  T_{-n} \le  T_{-1} \right) +   \mu\left(  T_{n} \le
  T_{-1} \right) \Bigr) 
\end{aligned}
\end{equation*}  
where we denote by $T_{\pm n}$ the first time for which $B_1 = \pm n /2$, and 
\ $T_{-1}$ the first time for which $ B_2 = -1$. 
It is well-known that the distributions of these hitting times are 
$$ dT_{\pm n}\left( \mu\right)  \: =  \: \frac{n}{4 \sqrt{2\pi t^3}}
\exp\left(\frac{-n^2}{8t}\right) {\bf 1}_{t > 0} \, dt, \quad 
  dT_{-1}\left( \mu \right) \: =  \: \frac{1}{ \sqrt{2\pi t^3}}
\exp\left(\frac{-1}{2t}\right) {\bf 1}_{t > 0} \, dt $$
A straightforward calculation provides 
\begin{equation*}
 \mu\left(  T_{\pm n} \le  T_{-1} \right) \:  = \: \int_{0 \le t_1 \le
   t_2} d T_{\pm n}\left( \mu \right)(t_1) \:\:  dT_{-1}\left( \mu
   \right)(t_2) \: \le \:  C (n^2+1)^{-1/2} 
\end{equation*}
which gives the result. 
\subsection{Optimality of the decay rate}
Theorem \ref{theosharp} shows that the boundary layer solution $v$
converges  {\em at least as} $y_2^{-1/2}$. One may
wonder if this result is optimal, that is if we can find
roughness distributions for which the speed of convergence is exactly given by
$y_2^{-1/2}$. In other words, is the constant $\sigma_{(0,0)}$ of the theorem
positive for some random distribution of roughness ? We have not so far 
 been able to show  optimality in this setting, but it can be established
 for the easier Dirichlet problem 
\begin{equation*}
\left\{
\begin{aligned}
\Delta u & = 0, \quad y_2 > 0 \\
       u & = \omega, \quad y_2 = 0 
\end{aligned}
\right.
\end{equation*}
where $\omega$ is a given boundary data. Although simpler, 
this system shares many features with the original system \eqref{BL}:
\begin{itemize}
\item If $\omega = \omega(y_1)$ is say $L$-periodic, 
the  solution $u(y)$  converges exponentially fast to the constant $\alpha :=
 L^{-1} \int_0^L \omega(y_1) dy_1$, as $y_2$ goes to infinity. 
\item If $\omega$ belongs as before to the probability space $(P,{\cal C},
 \pi)$, 
 one can show under assumption (H2) that 
$$ y_2  \, \E | u(\omega,0,y_2) - \alpha |^2 \: \xrightarrow[y_2
  \rightarrow +\infty]{} \: \sigma^2 \:  \ge   0,  \quad \alpha :=
  \E(\omega \mapsto \omega(0)).$$   
\end{itemize} 
Along the lines of \cite[pages 21-22]{Varadhan:2007}, 
we will exhibit a stationary measure $\pi$ for which $\sigma >
0$. Of course,  $\pi$ is  the law of the random process
$\varphi(\omega, y_1) := \omega(y_1)$, so that we just need to characterize
the random initial data.  Let $G(\omega,y_1)$  a gaussian 
random process, of zero mean and covariance $\rho(z_1-y_1)$, 
where $\rho \ge 0$  is a smooth even function with compact support.  Note
that such process exists: take $\rho = f * f$, with $f \ge 0$ an even
smooth function with compact support. Then, its Fourier transform satisfies 
$\hat{\rho} = |\hat{f}|^2 \ge 0$, which ensures the required positivity
property 
$$ \sum_{z_1,y_1} \, c(z_1) \,  c(y_1) \rho(z_1-y_1) \: =\ : \int_\R |
\sum_{z_1}  c(z_1) e^{i\xi z_1} |^2 \hat{\rho}(\xi) \, d\xi \: \ge \: 0. $$ 
for any family $c$ with compact support. Note moreover that this process
defines almost surely smooth functions of $y_1$: indeed, a simple
calculation yields  
$$ \E \left( \int_{[-R,R]} |\pa_{y_1}^k X(\cdot,y_1)|^2  \, dy_1 \right) 
\: = \: 2R \, (-1)^k \,  \rho^{(2k)}(0) \: < \: \infty $$
so that $X(\omega, \cdot)$ is almost surely in the space $H^k_{loc}(\R)$
and therefore smooth. Finally, we  introduce
$$ \varphi(\omega,y_1) \: = \:  F(X(\omega,y_1)) $$
for a smooth  increasing function $F$ with values in $(0,1)$. We stress that
$\varphi$ satisfies (H2) as $\rho$ has compact support. We will show that
the corresponding $\sigma$ is positive. Suppose {\em a contrario} that
$\sigma = 0 $. For $y_2 \ge 1$, we introduce the measure $\pi^{y_2}$
 associated to the gaussian process with variance $\rho(z_1-y_1)$ but mean 
$$m(z_1,y_2) \: = \: \int \rho(z_1 - y_1)  g(y_1,y_2) \, dy_1 $$
 where $g$ will be given later. Note that $\pi^{y_2}$ is associated to the
 random initial data 
$$ \varphi^{y_2}(\omega,y_1) \: := \:  F(X(\omega,y_1) + m(y_1,y_2)).$$
  Standard computation yields 
\begin{align*}
 R(\omega,y_2) \: := \: \frac{d\pi^{y_2}}{d \pi}(\omega) \: = & 
 \: \exp\Bigl(\int  g(y_1,y_2) \varphi(\omega,y_1) dy_1 \\
 & - \, \frac{1}{2}
 \,  \int\int \rho(z_1 - y_1) g(z_1,y_2) \, g(y_1,y_2) dz_1 \, dy_1
 \Bigr). 
\end{align*}
and 
$$ \int  |R(\omega,y_2)|^2 d\pi = e^{H(y_2)}, \quad H(y_2) \: := \: \int\int
\rho(z_1 - y_1) g(z_1,y_2) \, g(y_1,y_2) dz_1 \, dy_1. $$ 
If $\sigma=0$, then a simple Cauchy-Schwartz inequality 
$$ \int \sqrt{y_2} | u(\omega,0,y_2) - \alpha| d\pi^{y_2} \: \le \: 
 e^{\frac{1}{2} \, H(y_2)} \,  \int y_2 | u(\omega,0,y_2) - \alpha|^2 
d\pi $$
and {\em goes to zero as $y_2 \rightarrow +\infty$ if  $H$ is bounded from
 above}. 

\medskip
 Let $u$ and $v$ solutions associated to the initial data
 $\varphi(\omega,y_1)$ and $\varphi^{y_2}(\omega,y_1)$. As $m \ge 0$, by
 monotonicity, $v \ge u$. We can express $u$ and $v$ in terms of the Poisson
 Kernel,  so that 
$$ (v - u)(\omega,0,y_2) \: = \: \frac{2}{\pi} \int_\R \frac{y_2}{y_1^2 +
   y_2^2} \left( \varphi^{y_2}(\omega,y_1) - \varphi(\omega,y_1) \right)
dy_1. $$
Now, we define for $y_2 \ge 1$ 
$$ g(y_1,y_2) \: = \: \frac{1}{\sqrt{y_2}} \, G\left( \frac{y_1}{y_2}
\right), $$ 
where $G \ge 0$ has compact support, $G=1$ over $(-1,1)$. On one hand, 
with this definition of $g$, one can check that 
$$ \sup_{y_2 \ge 1} H(y_2) \: = \: \sup_{y_2 \ge 1} y_2^{-1} \int\int
\rho(z_1 - y_1)    G\left( \frac{z_1}{y_2}
\right) \,  G\left( \frac{y_1}{y_2}\right) dz_1 dy_1  \: < \: +\infty $$
On the other hand, one has 
\begin{align*}
& \int \sqrt{y_2} (v - u)(\omega,0,y_2) d\pi \:  \ge  \:  C  \int_\R \left( 
\int F'(X(\omega,0)) d\pi \right)  \frac{y_2}{y_1^2 + y_2^2} \, 
\sqrt{y_2} \, m(y_1,y_2) \, 
dy_1 \\
& \ge \: C' \,   \int_\R \left( \int F'(X(\omega,0)) d\pi \right) \, \left( 
\inf_{y_2 \ge 1} \inf_{|y_1| \le y_2} \int \rho(y_1
- s) G\left( \frac{s}{y_2}\right) ds \, \right) 
\left(  \int_{-y_2}^{y_2}   \frac{y_2}{y_1^2 + y_2^2} dy_1  \right) \\
& \ge \: C'' \, \int_\R \rho(s) ds > 0. 
\end{align*}
This implies that the quantity $ \int \sqrt{y_2} \,  (v(\omega,0,y_2) - \alpha)
\, d\pi $   does not go to zero, leading to a contradiction.

{\small
\def\cprime{$'$} \def\cprime{$'$} \def\cprime{$'$}

}

\section*{Acknowledgements}
The author warmly thanks S.R.S Varadhan for pointing to reference
\cite{Varadhan:2007},
 as well as Luis Silvestre and Thierry Levy for fruitful discussions.   
\section*{Appendix: Measurability of $\tilde{w}$}
We want to show here that 
$$ \tilde{w}(\omega,z) \: :=  \: \int_{y_2=0} G_\omega(z,y) \, e_1 \, dy, 
\quad  z  \mbox{ in } \Omega^{bl}(\omega),
  \quad \tilde{w}(\omega,z) := 0 \quad \mbox{ otherwise}  $$
defines a  measurable function from $P$ to $H^1_{loc}(\R^2)$. 
Let  $0 \le \varphi_n \le 1$ a sequence of 
 smooth functions with compact support,  $\: \varphi_n\vert_{(-n,n)} =
 1$. We define 
$$ w_n \: := \: \int _{y_2=0} G_\omega(z,y) \, (\varphi_n e_1) \, dy, 
\quad  z  \mbox{ in } \Omega^{bl}(\omega),
\quad w_n(\omega,z) := 0 \quad \mbox{ otherwise}.  $$
Note that $w_n$ is the (unique) solution of
\begin{equation} \label{BL3}
\left\{
\begin{aligned}
& -\Delta w_n + \na q_n = 0, 
\: x \in \Omega^{bl}\setminus \{ y_2 = 0\}, \\
& \div w_n = 0, \: x \in \Omega^{bl},\\
& w_n\vert_{\pa \Omega^{bl}} \: = \: 0,\\
& [ w_n ]\vert_{y_2 = 0} = 0, \quad [\pa_2 w_n - (0,q_n)]\vert_{y_2 = 0} = 
 (-\varphi_n,0),  
\end{aligned}
\right.
\end{equation}
satisfying $ \int_{\Omega^{bl}(\omega)} |\na w_n|^2 < +\infty$. 
By the dominated convergence theorem applied to the integral
formula, we get that  $w_n  \rightarrow \tilde{w}$
in $L^2_{loc}$. By the Cacciopoli inequality, the convergence is also
true in $H^1_{loc}$. Thus, we just have to show measurability of $w_n$. 

\medskip
Let us define
$$V \:  := \: \left\{ v \in \dot{H}^1(\R^2), \quad \div v = 0 \right\},
 \quad V_\omega \: := \:  \left\{ v
 \in V, \quad v\vert_{\R^2\setminus\Omega^{bl}(\omega)} = 0 \right\}.$$ 
 Following the lines of \cite[pages 15-16]{Basson:2007} 
it can be shown that the application 
$\omega \mapsto \pi(\omega)$, where $\pi(\omega) \in {\cal
  L}(V,V)$ is the orthogonal 
projection from $V$ to $V_\omega$, is measurable.  Now,
  $w_n$ is the unique fixed point of  the contraction 
$$ w  \mapsto \frac{1}{2}\pi(\omega) \left( w - v_n  \right) $$  
where $v_n$ is the unique function of $H^1(\R^2)$ satisfying 
$$ \int_{\R^2}  \na v_n \cdot  \na \phi \: = \: 6 \int_{y_2 = 0} \varphi_n \, 
\phi_1. $$
 The measurability of $w_n$ follows. 

\begin{thebibliography}{10}

\bibitem{Achdou:1998a}
{\sc Achdou, Y., Le~Tallec, P., Valentin, F., and Pironneau, O.}
\newblock Constructing wall laws with domain decomposition or asymptotic
  expansion techniques.
\newblock {\em Comput. Methods Appl. Mech. Engrg. 151}, 1-2 (1998), 215--232.
\newblock Symposium on Advances in Computational Mechanics, Vol.\ 3 (Austin,
  TX, 1997).

\bibitem{Achdou:1995}
{\sc Achdou, Y., Mohammadi, B., Pironneau, O., and Valentin, F.}
\newblock Domain decomposition \& wall laws.
\newblock In {\em Recent developments in domain decomposition methods and flow
  problems (Kyoto, 1996; Anacapri, 1996)}, vol.~11 of {\em GAKUTO Internat.
  Ser. Math. Sci. Appl.} Gakk\=otosho, Tokyo, 1998, pp.~1--14.

\bibitem{Achdou:1998}
{\sc Achdou, Y., Pironneau, O., and Valentin, F.}
\newblock Effective boundary conditions for laminar flows over periodic rough
  boundaries.
\newblock {\em J. Comput. Phys. 147}, 1 (1998), 187--218.

\bibitem{Amirat:2001a}
{\sc Amirat, Y., Bresch, D., Lemoine, J., and Simon, J.}
\newblock Effect of rugosity on a flow governed by stationary {N}avier-{S}tokes
  equations.
\newblock {\em Quart. Appl. Math. 59}, 4 (2001), 769--785.

\bibitem{Avellaneda:1987}
{\sc Avellaneda, M., and Lin, F.-H.}
\newblock Compactness methods in the theory of homogenization.
\newblock {\em Comm. Pure Appl. Math. 40}, 6 (1987), 803--847.

\bibitem{Avellaneda:1991}
{\sc Avellaneda, M., and Lin, F.-H.}
\newblock {$L\sp p$} bounds on singular integrals in homogenization.
\newblock {\em Comm. Pure Appl. Math. 44}, 8-9 (1991), 897--910.

\bibitem{Baladi:2001}
{\sc Baladi, V.}
\newblock Decay of correlations.
\newblock In {\em Smooth ergodic theory and its applications (Seattle, WA,
  1999)}, vol.~69 of {\em Proc. Sympos. Pure Math.} Amer. Math. Soc.,
  Providence, RI, 2001, pp.~297--325.

\bibitem{Basson:2007}
{\sc Basson, A., and G\'erard-Varet, D.}
\newblock Wall laws for fluid flows at a boundary with random roughness.
\newblock {\em Comm. Pure Applied Math., to appear} (2007). 
 Available at
\verb?http://www.dma.ens.fr/~dgerardv/publications.html?.


\bibitem{Bechert:1989}
{\sc Bechert, D., and Bartenwerfer, M.}
\newblock The viscous flow on surfaces with longitudinal ribs.
\newblock {\em J. Fluid Mech. 206}, 1 (1989), 105--129.

\bibitem{Bresch:2005}
{\sc Bresch, D., and G{\'e}rard-Varet, D.}
\newblock Roughness-induced effects on the quasi-geostrophic model.
\newblock {\em Comm. Math. Phys. 253}, 1 (2005), 81--119.

\bibitem{Bresch:2006}
{\sc Bresch, D., and Milisic, V.}
\newblock Higher order boundary layer correctors and wall laws derivation: a
  unified approach.
\newblock {\em Preprint arXiv:math/0611083 (2006)}.

\bibitem{Bouard:2007}
{\sc De~Bouard, A., Craig, W., Díaz-Espinosa, O., Guyenne, P., and C., S.}
\newblock Long wave expansions for water waves over random topography.
\newblock {\em Preprint arXiv:math/0506595} (2007).

\bibitem{Durrett:1996}
{\sc Durrett, R.}
\newblock {\em Probability: theory and examples}, second~ed.
\newblock Duxbury Press, Belmont, CA, 1996.

\bibitem{Galdi:1994}
{\sc Galdi, G.~P.}
\newblock {\em An introduction to the mathematical theory of the
  {N}avier-{S}tokes equations. {V}ol. {I}}, vol.~38 of {\em Springer Tracts in
  Natural Philosophy}.
\newblock Springer-Verlag, New York, 1994.
\newblock Linearized steady problems.

\bibitem{Gerard-Varet:2003b}
{\sc G{\'e}rard-Varet, D.}
\newblock Highly rotating fluids in rough domains.
\newblock {\em J. Math. Pures Appl. (9) 82}, 11 (2003), 1453--1498.

\bibitem{Giaquinta:1983}
{\sc Giaquinta, M.}
\newblock {\em Multiple integrals in the calculus of variations and nonlinear
  elliptic systems}, vol.~105 of {\em Annals of Mathematics Studies}.
\newblock Princeton University Press, Princeton, NJ, 1983.

\bibitem{Giaquinta:1982}
{\sc Giaquinta, M., and Modica, G.}
\newblock Nonlinear systems of the type of the stationary {N}avier-{S}tokes
  system.
\newblock {\em J. Reine Angew. Math. 330\/} (1982), 173--214.

\bibitem{Gilbarg:2001}
{\sc Gilbarg, D., and Trudinger, N.~S.}
\newblock {\em Elliptic partial differential equations of second order}.
\newblock Classics in Mathematics. Springer-Verlag, Berlin, 2001.
\newblock Reprint of the 1998 edition.

\bibitem{Jager:2001}
{\sc J{\"a}ger, W., and Mikeli{\'c}, A.}
\newblock On the roughness-induced effective boundary conditions for an
  incompressible viscous flow.
\newblock {\em J. Differential Equations 170}, 1 (2001), 96--122.

\bibitem{Jager:2003}
{\sc J{\"a}ger, W., and Mikeli{\'c}, A.}
\newblock Couette flows over a rough boundary and drag reduction.
\newblock {\em Comm. Math. Phys. 232}, 3 (2003), 429--455.

\bibitem{Jager:2000}
{\sc J{\"a}ger, W., Mikeli{\'c}, A., and Neuss, N.}
\newblock Asymptotic analysis of the laminar viscous flow over a porous bed.
\newblock {\em SIAM J. Sci. Comput. 22}, 6 (2000), 2006--2028 (electronic).

\bibitem{Jikov:1994}
{\sc Jikov, V.~V., Kozlov, S.~M., and Ole{\u\i}nik, O.~A.}
\newblock {\em Homogenization of differential operators and integral
  functionals}.
\newblock Springer-Verlag, Berlin, 1994.
\newblock Translated from the Russian by G. A. Yosifian [G. A. Iosif\cprime
  yan].

\bibitem{Luchini:1995}
{\sc Luchini, P.}
\newblock Asymptotic analysis of laminar boundary-layer flow over finely
  grooved surfaces.
\newblock {\em European J. Mech. B Fluids 14}, 2 (1995), 169--195.

\bibitem{Varadhan:2007}
{\sc Varadhan, S., and Zygouras, N.}
\newblock Behavior of the solution of a random semilinear heat equation.
\newblock {\em Comm. Pure Applied Math., to appear} (2007). 
 Available at
  \verb?http://www-rcf.usc.edu/~zygouras/research.html?.

\bibitem{Varadhan:1968}
{\sc Varadhan, S. R.~S.}
\newblock {\em Stochastic processes}.
\newblock Notes based on a course given at New York University during the year
  1967/68. Courant Institute of Mathematical Sciences New York University, New
  York, 1968.

\end{thebibliography}
\end{document}